\documentclass{amsart}
\usepackage[utf8]{inputenc}
\usepackage{amsthm,amsmath,amssymb,tikz-cd,mathrsfs,enumitem}

\makeatletter
\@namedef{subjclassname@2010}{%
  \textup{2010} Mathematics Subject Classification}
\makeatother

\title[$C^*$-modules over $\Sigma^*$-algebras]{Hilbert $C^*$-modules over $\Sigma^*$-algebras}
\author[C. A. Bearden]{Clifford A. Bearden}
\address{Department of Mathematics\\ University of Houston\\ Houston, TX 77204-3008, USA}
\email{cabearde@math.uh.edu}

\date{}

\newtheorem{thm}{Theorem}[section]
\newtheorem{prop}[thm]{Proposition}
\newtheorem{lemma}[thm]{Lemma}
\newtheorem{cor}[thm]{Corollary}

\theoremstyle{definition}
\newtheorem{defn}[thm]{Definition}

\newtheorem{note}[thm]{Note}

\newtheorem{example}[thm]{Example}
\newtheorem{remark}[thm]{Remark}

\newcommand{\C}{\mathbb{C}}
\newcommand{\N}{\mathbb{N}}
\newcommand{\R}{\mathbb{R}}
\newcommand{\K}{\mathbb{K}}
\newcommand{\B}{\mathbb{B}}
\newcommand{\Hil}{\mathcal{H}}
\newcommand{\Kil}{\mathcal{K}}
\newcommand{\bw}{\Sigma^*}
\newcommand{\swb}{\Sigma^*_{\mathfrak{B}}\text{-}}
\newcommand{\swa}{\Sigma^*_A\text{-}}
\newcommand{\Ball}{\text{Ball}}
\newcommand{\Ker}{\text{Ker}}
\newcommand{\id}{\text{id}}

\newcommand{\fA}{\mathfrak{A}}
\newcommand{\fB}{\mathfrak{B}}
\newcommand{\fC}{\mathfrak{C}}
\newcommand{\fD}{\mathfrak{D}}
\newcommand{\fX}{\mathfrak{X}}
\newcommand{\fY}{\mathfrak{Y}}
\newcommand{\fZ}{\mathfrak{Z}}
\newcommand{\Bw}{\mathscr{B}}

\newcommand{\xp}{\mathscr{B}(X)}

\begin{document}

\renewcommand{\thefootnote}{}
%
%

\subjclass[2010]{Primary 46L08; Secondary 28A20}

\keywords{$C^*$-module, $W^*$-module, $\Sigma^*$-algebra, TRO, corner, selfdual, countably generated}

\begin{abstract}
A \emph{$\Sigma^*$-algebra} is a concrete $C^*$-algebra that is sequentially closed in the weak operator topology. We study an appropriate class of $C^*$-modules over $\Sigma^*$-algebras analogous to the class of $W^*$-modules (selfdual $C^*$-modules over $W^*$-algebras), and we are able to obtain $\Sigma^*$-versions of virtually all the results in the basic theory of $C^*$- and $W^*$-modules. In the second half of the paper, we study modules possessing a weak sequential form of the condition of being countably generated. A particular highlight of the paper is the ``$\Sigma^*$-module completion," a $\Sigma^*$-analogue of the selfdual completion of a $C^*$-module over a $W^*$-algebra, which has an elegant uniqueness condition in the countably generated case.
\end{abstract}

\maketitle

\section{Introduction}

Hilbert $C^*$-modules (also called Hilbert modules, and which we simply call $C^*$-modules) are simultaneous generalizations of $C^*$-algebras, Hilbert spaces, and certain types of vector bundles. They are an amazingly versatile tool used in a broad range of subfields of operator algebra theory---for example, the theory of Morita equivalence, Kasparov's $KK$-theory and its applications in noncommutative geometry, quantum group theory, and operator space theory (see \cite[Section 8.6]{BLM} for the latter).

An important subclass of $C^*$-modules is the class of selfdual $C^*$-modules (see Definition \ref{sddef}) over $W^*$-algebras, i.e.\ the \emph{$W^*$-modules}. Historically, $W^*$-modules (introduced in 1973 by Paschke in \cite{Pas}) were the first $C^*$-modules to garner wide-reaching attention, but today they seem less well-known and perhaps under-exploited. Compared to the general theory of $C^*$-modules, the theory of $W^*$-modules is much more elegant and similar to that of Hilbert spaces, in large part due to powerful ``orthogonality" properties automatically present in $W^*$-modules.

Between the classes of $C^*$-algebras and $W^*$-algebras is the class of \emph{$\bw$-algebras}---first studied by Davies in \cite{Dav68}, these are defined as concrete $C^*$-algebras that are sequentially closed in the weak operator topology (abbreviated WOT from here on). It is the purpose of this paper to explore the ``appropriate" class of $C^*$-modules over $\bw$-algebras in analogy with the way that $W^*$-modules are the ``appropriate" class of $C^*$-modules over $W^*$-algebras.

This paper is broken up into three sections. In the first, we quickly survey some background facts about $C^*$-modules, $W^*$-modules, and $\bw$-algebras.

In the second, we define our class of ``$\bw$-modules" and prove general results analogous to many of the basic results in $C^*$- and $W^*$-module theory. In particular, we show that $\bw$-modules correspond with the ternary rings of operator (TROs) that are WOT sequentially closed and with corners of $\bw$-algebras in the same way that $C^*$-modules (resp.\ $W^*$-modules) correspond with norm-closed (resp.\ weak*-closed) TROs and with corners of $C^*$-algebras (resp.\ $W^*$-algebras). The other main highlight of this section is the ``$\bw$-module completion" of a $C^*$-module over a $\bw$-algebra, in analogy with the selfdual completion of a $C^*$-module over a $W^*$-algebra.

In the final section, we study the subclass of ``$\swb$countably generated" $\bw$-modules, and are able to prove many satisfying results about these---for example, that all $\swb$countably generated $\bw$-modules are selfdual. As expected, there is also a weak sequential version of Kasparov's stabilization theorem that holds in this case.

This work was inspired in large part by Hamana's paper \cite{Hamana}, in which he studies selfdual $C^*$-modules over monotone complete $C^*$-algebras. Theorems 1.2, 2.2, and 3.3 of that paper indicate that selfdual $C^*$-modules are the ``appropriate" class of $C^*$-modules over monotone complete $C^*$-algebras (and his interesting ``conversely" statement in Theorem 2.2 seems to indicate that monotone complete $C^*$-algebras are the ``appropriate" class of coefficient $C^*$-algebras over which to consider selfdual $C^*$-modules). We do not have in the case of $\bw$-modules the existence of an ``orthonormal basis," which is Hamana's main technical tool in \cite{Hamana}, so most of his proof techniques do not work for us, but the overarching philosophy of what we have tried to accomplish is very much in line with that of Hamana's work. For more work on the subject of $C^*$-modules over monotone complete $C^*$-algebras, see the paper \cite{Frank} by M.\ Frank.

Also somewhat related to the present work is the noncommutative semicontinuity theory initiated by Akemann and Pedersen in \cite{AP} and developed further by Brown in \cite{BrownSem} (see also \cite{BrownDS,BrownLarge}). Though the present work has seemingly little to do with this theory (we do not deal with monotone limits, and in this work, the universal representation is mentioned only as an example setting), we were first drawn into this investigation by Brown's mentioning in \cite{BrownLarge} that the monotone sequentially closed $C^*$-algebra generated by the set of semicontinuous elements in the second dual of a $C^*$-algebra is seemingly the most natural noncommutative analogue of the space of bounded Borel functions on a locally compact Hausdorff space. See Note \ref{semiconts qs} for a short discussion of some related interesting open problems.

Results in this paper were announced at the Workshop in Analysis and Probability at Texas A\&M University in July 2015 and the mini-workshop ``Operator Spaces and Noncommutative Geometry in Interaction" at Oberwolfach in February 2016. The author would like to thank the organizers of these conferences, as well as the NSF for funding during part of the time research was carried out for this paper. We would also like to thank Larry Brown for helpful discussions during a visit of his to Houston in September 2014.

This work constitutes a significant part of the author's doctoral thesis at the University of Houston. The author would like to express his deep gratitude for the guidance of his Ph.D. advisor, David Blecher, who offered many suggestions and helped answer a couple of questions in preliminary versions of this paper, and without whose guidance and support none of this would have been possible.

\section{Background}

In this section we fix our notation and review the basic definitions and results for $C^*$-modules, $W^*$-modules, and $\bw$-algebras. Since the basic theory of $C^*$-modules is well-known and covered in many texts, we will be brief here. We generally refer to \cite[Chapter 8]{BLM} for notation and results; other references include \cite{Lan}, \cite{RW}, \cite[Chapter 15]{WO}, and \cite[Section II.7]{Bl}.

Loosely speaking, a (right or left) module $X$ over a $C^*$-algebra $A$ is called a \emph{(right or left) $C^*$-module over $A$} if it is equipped with an ``$A$-valued inner product" $\langle \cdot | \cdot \rangle : X \times X \to A$ and is complete in the canonical norm induced by this inner product. If $X$ is a right $C^*$-module, the inner product is taken to be linear and $A$-linear in the second variable and conjugate-linear in the first variable, and vice versa for left modules. When unspecified, ``$C^*$-module" should be taken to mean ``right $C^*$-module."

If $X$ and $Y$ are two $C^*$-modules over $A$, $B_A(X,Y)$ denotes the Banach space of bounded $A$-module maps from $X$ to $Y$ with operator norm; $\mathbb{B}_A(X,Y)$ denotes the closed subspace of adjointable operators; and $\mathbb{K}_A(X,Y)$ denotes the closed subspace generated by operators of the form $|y \rangle \langle x| := y \langle x | \cdot \rangle$ for $y \in Y,$ $x \in X.$ If $X=Y$, the latter two of these spaces are $C^*$-algebras, and in this case, $X$ is a left $C^*$-module over $\K_A(X)$ with inner product $|\cdot \rangle \langle \cdot|.$

In this paper, we will be concerned with modules over $\bw$-algebras, a class of $C^*$-algebras with an extra bit of structure that may be viewed abstractly, but is most easily captured by fixing a faithful representation of a certain type on a Hilbert space. Reflecting this view, $C^*$-modules over $\bw$-algebras are also most easily studied when viewed under a representation induced by a fixed representation of the coefficient $\bw$-algebra. There is a well-known general procedure for taking a representation of the coefficient $C^*$-algebra of a $C^*$-module and inducing a representation of the $C^*$-module and many of the associated mapping spaces mentioned in the previous paragraph. The following paragraph and proposition describe this construction and its relevant features.

If $A$ is a $C^*$-algebra represented nondegenerately on a Hilbert space $\Hil$ and $X$ is a $C^*$-module over $A$, we may consider $\Hil$ as a left module over $A$ and take the algebraic module tensor product $X \odot_A \Hil.$ This vector space admits an inner product determined by the formula $\langle x \otimes \zeta, y \otimes \eta \rangle = \langle \zeta, \langle x|y \rangle \eta \rangle$ for simple tensors (see \cite[Proposition 4.5]{Lan} for details), and we may complete $X \odot_A \Hil$ in the induced norm to yield a Hilbert space $X \otimes_A \Hil.$ Considering $A$ as a $C^*$-module over itself and taking the $C^*$-module direct sum $X \oplus A$, there is a canonical corner-preserving embedding of $B_A(X \oplus A)$ into $B((X \otimes_A \Hil) \oplus^2 \Hil)$ which allows us to concretely identify many of the associated spaces of operators between $X$ and $A$ with spaces of Hilbert space operators between $\Hil$ and $X \otimes_A \Hil$---this is the content of the following proposition. All of the pieces of this proposition can be found in the textbooks mentioned above.

Recall that for a nondegenerately-acting $C^*$-algebra $A \subseteq B(\Hil)$, the multiplier algebra $M(A)$ may be identified with the space $\{T \in B(\Hil) : TA \subseteq A \text{ and } AT \subseteq A\},$ and the left multiplier algebra $LM(A)$ with $\{T \in B(\Hil): TA \subseteq A\}.$ For a right $C^*$-module $X$ over $A$, we write $``\overline{X}"$ to denote the adjoint $C^*$-module (see \cite[8.1.1]{BLM})---this is a left $C^*$-module over $A$.

\begin{prop} \label{bigdiagram}
If $X$ is a $C^*$-module over a nondegenerate $C^*$-algebra $A \subseteq B(\Hil)$, then there are canonical maps making the following diagram commute:

$$
\begin{tikzcd}
\K_A(X \oplus A) \arrow{r} \arrow{d}
& \left[ \begin{matrix} \K_A(X) & \K_A(A,X) \\ \K_A(X,A) & \K_A(A) \end{matrix} \right] \arrow{r} \dar[shorten <= 1ex, shorten >= 1ex]
& \left[ \begin{matrix} \K_A(X) & X \\ \overline{X} & A \end{matrix} \right] \dar[shorten <= 1ex, shorten >= 1ex]
\\
\B_A(X \oplus A) \arrow{r} \arrow{d}
& \left[ \begin{matrix} \B_A(X) & \B_A(A,X) \\ \B_A(X,A) & \B_A(A) \end{matrix} \right] \arrow{r} \dar[shorten <= 1ex, shorten >= 1ex]
& \left[ \begin{matrix} M(\K_A(X)) & * \\ * & M(A) \end{matrix} \right] \dar[shorten <= 1ex, shorten >= 1ex] 
\\
B_A(X \oplus A) \arrow{r} \arrow{d}
& \left[ \begin{matrix} B_A(X) & B_A(A,X) \\ B_A(X,A) & B_A(A) \end{matrix} \right] \arrow{r} \dar[shorten <= 1ex, shorten >= 1ex]
& \left[ \begin{matrix} LM(\K_A(X)) & * \\ * & LM(A) \end{matrix} \right] \dar[shorten <= 1ex, shorten >= 1ex] 
\\
B((X \otimes_A \Hil) \oplus^2 \Hil) \arrow{r} & \left[ \begin{matrix} B(X \otimes_A \Hil) & B(\Hil, X \otimes_A \Hil) \\ B(X \otimes_A \Hil,\Hil) & B(\Hil) \end{matrix} \right] \arrow{r}
& \left[ \begin{matrix} B(X \otimes_A \Hil) & B(\Hil, X \otimes_A \Hil) \\ B(X \otimes_A \Hil,\Hil) & B(\Hil) \end{matrix} \right]
\end{tikzcd}
$$

The horizontal maps in the third row are Banach algebra isomorphisms, and all other horizontal maps are $*$-isomorphisms. All vertical maps are isometric homomorphisms, and in the diagram with the third row deleted, all vertical maps are isometric $*$-homomorphisms. In the diagram with the first column deleted, all maps are corner-preserving.

Furthermore, if we define $\mathcal{L}(X):=\K_A(X \oplus A)$ (called the \emph{linking $C^*$-algebra of $X$}), then $M(\mathcal{L}(X))=\B_A(X \oplus A)$ and $LM(\mathcal{L}(X))=B_A(X \oplus A).$
\end{prop}

\begin{note}
We will often use the proposition above many times in the sequel, often without mention and often without distinguishing between a $C^*$-module operator and its image as a Hilbert space operator. That said, we will sometimes have two $C^*$-algebras $\fA \subseteq B(\Kil)$ and $\fB \subseteq B(\Hil)$ and a bimodule $X$ that is a left $C^*$-module over $\fA$ and a right $C^*$-module over $\fB;$ in this case, it is important to distinguish whether we are viewing $X$ as embedded in $B(\Hil, X \otimes_\fB \Hil)$ or in $B(\overline{X} \otimes_\fA \Kil, \Kil)$ (see Note \ref{dist notation}).
\end{note}

\begin{defn} \label{sddef}
A right $C^*$-module $X$ over a $A$ is called \emph{selfdual} if every bounded $A$-module map $X \to A$ is of the form $\langle x|\cdot \rangle$ for some $x \in X.$ A \emph{$W^*$-module} is a selfdual $C^*$-module over a $W^*$-algebra.
\end{defn}

There are many beautiful characterizations of $W^*$-modules among $C^*$-modules. Most elegantly, a $C^*$-module over a $W^*$-algebra is a $W^*$-module if and only if it has a Banach space predual (this was originally proved in \cite{Zettl} and \cite{EOR}, or see \cite[Corollary 3.5]{BM} for another proof). For the purposes of this paper, the following characterization may be taken as motivation:

\begin{prop} \label{w* mod}
A $C^*$-module $Y$ over a von Neumann algebra $M \subseteq B(\Hil)$ is a $W^*$-module if and only if the canonical image of $Y$ in $B(\Hil, Y \otimes_M \Hil)$ is weak*-closed.
\end{prop}

\proof ($\Rightarrow$) Assume $Y$ is a $W^*$-module. By the Krein-Smulian theorem, it suffices to prove that if $(y_\lambda)$ is a bounded net in $Y$ such that $y_\lambda \xrightarrow{w^*} T$ in $B(\Hil,X\otimes_M \Hil),$ then $T \in Y$. If we have such a net $(y_\lambda)$ and operator $T$, then for any $x \in Y,$ $\langle y_\lambda \otimes \zeta, x \otimes \eta \rangle = \langle \zeta, \langle y_\lambda|x \rangle \eta \rangle$ is convergent for all $\zeta, \eta \in \Hil,$ hence $\langle y_\lambda|x \rangle$ converges WOT to some $a_x \in M$, and since the WOT and weak* topology on $B(\Hil)$ coincide on bounded sets, we have $\langle y_\lambda|x \rangle \xrightarrow{w^*} a_x$. Since $Y$ is a $W^*$-module, the map $x \mapsto a_x$ has the form $\langle y|\cdot \rangle$ for some $y \in Y,$ and it follows easily that $T=y$ in $B(\Hil, Y \otimes_M \Hil).$

($\Leftarrow$) Assume the latter condition, let $\varphi \in B_M(Y,M),$ and let $(e_t)$ be a cai (contractive approximate identity) for $\K_M(Y).$ For each $t$, $\varphi e_t \in \K_M(Y,M)$, and so there is a $y_t \in Y$ such that $\varphi e_t = \langle y_t | \cdot \rangle$ (by the top right isomorphism in Proposition \ref{bigdiagram}). By assumption, $Y$ is a dual Banach space, and so $(y_t)$ has a weak*-convergent subnet $y_{t_s} \xrightarrow{w^*} y.$

Using Cohen's factorization theorem (\cite[A.6.2]{BLM}) to write any $x \in Y$ as $x=Kx'$ for $K \in \K_M(Y)$ and $x' \in Y,$ we have $e_t(x)= e_t(Kx') = (e_t K)(x') \to Kx' = x$ in norm in $Y$, so $(\varphi e_t)(x) = \varphi(e_t x) \to \varphi(x)$ in norm in $M.$ Hence $(\varphi e_t)(x \otimes \zeta) = (\varphi e_t)(x)(\zeta) \to \varphi(x)(\zeta)=\varphi(x \otimes \zeta)$ in $\Hil$ for all $x \in Y$ and $\zeta \in \Hil$. Since $(\varphi e_t)$ is bounded and the simple tensors are total in $Y \otimes_M \Hil$, a triangle inequality argument shows that $(\varphi e_t)$ converges in the SOT (strong operator topology), hence weak*, to $\varphi$ in $B(Y \otimes_M \Hil,\Hil)$.

Since $y_{t_s} \xrightarrow{w^*} y$ in $B(\Hil, Y \otimes_M \Hil)$, we have $$\langle \langle y_{t_s}|\cdot \rangle(x \otimes \zeta), \eta \rangle = \langle x \otimes \zeta, y_{t_s}(\eta) \rangle \to \langle x \otimes \zeta, y(\eta) \rangle = \langle \langle y|\cdot \rangle(x \otimes \zeta), \eta \rangle.$$ Since $(\langle y_{t_s}| \cdot \rangle)$ is bounded, another triangle inequality argument as in the previous paragraph gives $\langle y_{t_s}| \cdot \rangle \xrightarrow{WOT} \langle y| \cdot \rangle,$ so that $\langle y_{t_s}| \cdot \rangle \xrightarrow{w^*} \langle y| \cdot \rangle$ by boundedness again.

Since $\varphi e_t = \langle y_t | \cdot \rangle$, we may combine the previous two paragraphs to conclude that $\varphi = \langle y| \cdot \rangle.$
\endproof



\begin{defn}
A \emph{(concrete) $\bw$-algebra} is a nondegenerate $C^*$-algebra $\fB \subseteq B(\Hil)$ that is closed under limits of WOT-convergent sequences, i.e.\ whenever $(b_n)$ is a sequence in $\fB$ that converges in the weak operator topology of $B(\Hil)$ to an operator $T$, then $T \in \fB.$

For $\bw$-algebras $\fA \subseteq B(\Kil), \fB \subseteq B(\Hil)$, a $*$-homomorphism $\varphi : \fA \to \fB$ is called a \emph{$\bw$-homomorphism} if $a_n \xrightarrow{WOT} a$ in $\fA$ implies $\varphi(a_n) \xrightarrow{WOT} \varphi(a)$ in $\fB.$ If additionally $\varphi$ is a $*$-isomorphism and $\varphi^{-1}$ is a $\bw$-homomorphism, then $\varphi$ is called a \emph{$\bw$-isomorphism}.
\end{defn}

The following simple observation is well-known. We do not directly apply this result in the present work, but it provides an alternate definition for $\bw$-algebras and seems noteworthy.

\begin{lemma} \label{WOTsc=wsc}
A sequence in $B(\Hil)$ is weak*-convergent if and only if it is WOT-convergent. Hence a $C^*$-algebra $A \subseteq B(\Hil)$ is a $\bw$-algebra if and only if it is sequentially closed in the weak* topology of $B(\Hil)$.
\end{lemma}

\proof
One direction is obvious. The other follows by applying the uniform boundedness principle (twice) to see that a WOT-convergent sequence is automatically bounded, hence also weak*-convergent.
\endproof

One may also discuss abstract $\bw$-algebras, i.e.\ $C^*$-algebras that admit a faithful representation as a concrete $\bw$-algebra. In his original paper on $\bw$-algebras (\cite{Dav68}), E. B. Davies proved a characterization theorem for when an abstract $C^*$-algebra $A$ equipped with a collection of pairs $((a_n),a)$, each consisting of a sequence $(a_n) \subseteq A$ and an element $a \in A$, admits a faithful representation as a $\bw$-algebra whose WOT-convergent sequences and their limits are prescribed by the collection $\{((a_n),a)\}.$ These abstract $\bw$-algebras may also be described by replacing the collection $\{((a_n),a)\}$ with an appropriate subspace of the dual space $A^*$, or an appropriate closed convex subset of the state space $S(A)$. The latter perspective was mentioned by Davies in the original paper and is the underlying point of view in Dang's paper \cite{Dang}. More explicitly, Dang defines a $\bw$-algebra to be a pair $(A,S)$, where $A$ is a $C^*$-algebra and $S$ is a subset of $S(A)$ such that:
\begin{enumerate}
\item[\rm{(1)}] if $\varphi \in S$ and $a \in A$ with $\varphi(a^*a)=1,$ then $\varphi(a^* \cdot a) \in S;$
\item[\rm{(2)}]  if $\psi$ is a state on $A$ such that $\psi(a_n)$ converges for all sequences $(a_n)$ in ${}^\sigma S = \{(a_n) \in \ell^\infty(A) : \varphi(a_n) \text{ converges for all } \varphi \in S\},$ then $\psi \in S;$
\item[\rm{(3)}]  if $a \in A$ is nonzero, then $\varphi(a) \not= 0$ for some $\varphi \in S$;
\item[\rm{(4)}]  if $(a_n) \in {}^\sigma S,$ then there is an $a \in A$ such that $\varphi(a_n) \to \varphi(a)$ for all $\varphi \in S.$
\end{enumerate}
Elementary operator theoretic arguments show that if $A$ is WOT sequentially closed, then the collection of WOT sequentially continuous states meets these requirements. Conversely, as Dang points out, one may use a slight modification of the polarization identity ($b^* x a = \frac 14 \sum_{k=0}^3 i^k (a+i^k b)^* x (a + i^k b)$ for $a,x,b \in A$), to check that if $(A,S)$ is a $\bw$-algebra in Dang's sense, then $(A,{}^\sigma S)$ is a $\bw$-algebra in Davies' sense, so that by Davies' result, $A$ admits a representation as a $\bw$-algebra in our sense.

A similar class of $C^*$-algebras was studied by Pedersen in several papers (see \cite[Section 4.5]{Ped} for the main part of the theory and more references). He studied ``Borel $*$-algebras," which are concrete $C^*$-algebras closed under limits of bounded monotone sequences of selfadjoint elements. In some ways, Borel $*$-algebras are more technically forbidding (e.g.\ compare Proposition \ref{closure is $C^*$} below to \cite[Theorem 4.5.4]{Ped}), but in other ways they seem nicer---for example, it seems to be an open question whether or not a $*$-isomorphism between $\bw$-algebras is always a $\bw$-isomorphism, but it is easy to see that the analogous statement for Borel $*$-algebras is true.

For any subset $S \subseteq B(\Hil)$, denote by $\Bw(S)$ the smallest WOT sequentially closed subset of $B(\Hil)$ containing $S$. Such a set exists since the intersection of any two WOT sequentially closed subsets is also WOT sequentially closed. If there is ambiguity (for example if we represent a $C^*$-algebra on two different Hilbert spaces), we add a subscript: $\Bw_\Hil(S).$ These closures provide many examples of $\bw$-algebras:

\begin{prop} \label{closure is $C^*$}
If $A \subseteq B(\Hil)$ is a nondegenerate $C^*$-algebra, then $\Bw(A)$ is a $\bw$-algebra.
\end{prop}

\proof
(\cite[Lemma 2.1]{Dav68}) Fix $a \in A,$ and let $S = \{b \in \Bw(A): ab \in \Bw(A)\}.$ Clearly $S$ is WOT sequentially closed and contains $A$, so $S=\Bw(A).$ Hence $ab \in \Bw(A)$ for all $a \in A$ and $b \in \Bw(A).$ Similar tricks show that $bc \in \Bw(A)$ for all $b,c \in \Bw(A)$ and that $\Bw(A)$ is a $*$-invariant subspace of $B(\Hil).$ Since $\Bw(A)$ is also evidently norm-closed, the result follows.
\endproof

\begin{example} \label{bunch of examps}
\begin{enumerate}
\item[(1)]  Every von Neumann algebra is clearly a $\bw$-algebra. Conversely, if $\Hil$ is separable, then it follows from Pedersen's up-down theorem (\cite[2.4.3]{Ped}) that every Borel $*$-algebra, hence every $\bw$-algebra, in $B(\Hil)$ is a von Neumann algebra. (Kadison first proved this fact for $\bw$-algebras in an appendix to \cite{Dav68}.)

\item[(2)]  If $\Hil$ is a Hilbert space, then the ideal $\mathscr{S}$ of operators in $B(\Hil)$ with separable range is the $\bw$-algebra $\Bw(\K(\Hil))$, which is of course not unital if $\Hil$ is not separable. Indeed, it is a short exercise to see that every operator with separable range is a SOT-limit of a sequence of finite rank operators. Conversely, by basic operator theory, every compact operator has separable range. To see that $\mathscr{S}$ is WOT sequentially closed, suppose $(T_n)$ is a sequence in $\mathscr{S}$ converging in the WOT to $T \in B(\Hil).$ Then $P:= \vee_n r(T_n)$ (where $r(T_n)$ denotes the projection onto $\overline{\text{Ran}(T_n)}$) is a projection with separable range, and $T_n=PT_n \xrightarrow{WOT} PT.$ Hence $PT=T,$ so $T$ has separable range.

\item[(3)]  Let $A$ be a $C^*$-algebra considered as a concrete $C^*$-algebra in its universal representation $A \subseteq B(\Hil_u)$. The $\bw$-algebra $\Sigma^*(A):=\Bw(A)$ obtained here is called the \emph{Davies-Baire envelope of $A$} (following the terminology of \cite{SaitoWright}). It was proved by Davies in \cite[Theorem 3.2]{Dav68} that $\Sigma^*(A)$ is $\bw$-isomorphic to $\Bw_{\Hil_a}(A)$, where $A \hookrightarrow B(\Hil_a)$ is the atomic representation of $A$. 

\item[(4)]  (\cite[Corollary 3.3]{Dav68}, \cite[4.5.14]{Ped}) Let $X$ be a locally compact Hausdorff space. By basic $C^*$-algebra theory, the atomic representation of the commutative $C^*$-algebra $C_0(X)$ is the embedding of $C_0(X)$ into $B(\ell^2(X))$ as multiplication operators. By the last statement in the previous example, $\Sigma^*(C_0(X))$ may be identified with the WOT sequential closure of $C_0(X)$ in $B(\ell^2(X))$. This closure is evidently contained in the copy of the space of all bounded functions on $X$, $\ell^\infty(X)$, in $B(\ell^2(X))$. Since WOT-convergence of sequences in $\ell^\infty(X) \subseteq B(\ell^2(X))$ coincides with pointwise convergence of bounded sequences of functions, we may identify $\Sigma^*(C_0(X))$ with the space of functions known classically (sometimes) as the bounded Baire functions on $X$ (in the sense of \cite[6.2.10]{Ped NOW} or \cite{Halmos}). Recall two well-known classical facts about the Baire functions: (1) if $X$ is second countable, the Baire functions and the Borel-measurable functions on $X$ coincide, and (2) $X$ is $\sigma$-compact if and only if the constant functions are Baire. Thus $\Sigma^*(C_0(X))$ for non-$\sigma$-compact $X$ provides another example of a nonunital $\bw$-algebra.


\item[(5)]  If $A$ is a separable $C^*$-algebra and $\phi$ is a faithful state on $A$, then the GNS construction gives a faithful representation of $A$ as operators on a separable Hilbert space $\Hil_\phi$. By (1) above, $\Bw_{\Hil_\phi}(A)$ is the weak*-closure of $A$ in $B(\Hil_\phi).$ In particular, if $A=C(X)$ for a second countable compact Hausdorff space $X$, and $\mu$ is a finite positive Borel measure on $X$ such that $\int f d \mu >0$ for all nonzero positive $f \in C(X)$ (e.g.\ take $\mu$ to be Lebesgue measure on $X=[0,1]$), then by basic measure theory (see e.g.\ \cite[Example 4.1.2]{Murphy}), $\Bw_{L^2(X,\mu)}(C(X)) = L^\infty(X,\mu).$
\end{enumerate}
\end{example}

\begin{note}[Open Questions] \label{semiconts qs}
(Cf.\ \cite{Ped} 4.5.14, \cite[Section 5.3.1]{SaitoWright}) Though we do not address these in the present work, there are some interesting and natural open questions about $\bw$-algebras and similar classes of $C^*$-algebras.

As mentioned above Proposition \ref{closure is $C^*$}, it appears to be unknown whether or not every $*$-isomorphism between $\bw$-algebras is a $\bw$-isomorphism (or even if $*$-isomorphic $\bw$-algebras are necessarily $\bw$-isomorphic).

Related to (3) in Example \ref{bunch of examps}, if $A  \subseteq B(\Hil_u)$ is a $C^*$-algebra in its universal representation, it is  unknown whether or not one must have $\Bw(A)=\Bw^m(A)$, where the latter refers to the monotone sequential closure of $A$ (that is, $\Bw^m(A) = \Bw^m(A_{sa}) + i \Bw^m(A_{sa}),$ where $\Bw^m(A_{sa})$ is the smallest subset of $B(\Hil_u)_{sa}$ containing $A_{sa}$ and closed under limits of bounded increasing sequences). Clearly the inclusion $\Bw^m(A) \subseteq \Bw(A)$ always holds. Pedersen proved that $\Bw^m(A)$ is always a $C^*$-algebra (\cite[4.5.4]{Ped}) and that the equation $\Bw(A)=\Bw^m(A)$ does hold if $A$ is type I (\cite[Section 6.3]{Ped}).

One may also replace the monotone sequential closure $\Bw^m(A)$ in the paragraph above with a number of variants---for example, the SOT sequential closure of $A$, $\Bw^s(A)$. Clearly $\Bw^s(A)$ lies between $\Bw^m(A)$ and $\Bw(A)$, but as far as we know, the questions of whether or not $\Bw^s(A)$ always equals $\Bw(A)$ or $\Bw^m(A)$ are still open. (Note that by Lemma \ref{WOTsc=wsc} the weak* sequential closure of $A$ coincides with $\Bw(A).$)

In fact, as far as we can tell, there is no known example of \emph{any} monotone sequentially closed $C^*$-algebra that is not WOT sequentially closed (or a WOT sequentially closed $C^*$-algebra that is not SOT sequentially closed).

Somewhat similar in spirit is the interesting open question of whether or not $A_{sa}^m$ (the set of limits in $A_{sa}^{**}$ of bounded increasing nets in $A_{sa}$) is always norm-closed. Brown proved in \cite{BrownSem} (Corollary 3.25) that this does hold if $A$ is separable. See \cite{BrownDS} for an insightful discussion on this problem.
\end{note}

We now briefly record a few basic facts about $\bw$-algebras that we will use later.

\begin{prop} \label{polar decomp}
Let $T$ be an operator in a $\bw$-algebra $\fB \subseteq B(\Hil).$ If $T=U|T|$ is the polar decomposition of $T$, then $U \in \fB.$
\end{prop}

\proof
See \cite[Lemma 2.1]{Dav69} or \cite[4.5.16]{Ped}.
\endproof

\begin{prop} \label{sp proj}
If $\fB \subseteq B(\Hil)$ is a unital $\bw$-algebra and $x$ is a selfadjoint element in $\fB,$ then $f(x) \in \fB$ for all bounded Borel functions $f: \R \to \C.$
\end{prop}

\proof
This may be proved by a mild modification of \cite[4.5.7]{Ped}.
\endproof

\begin{prop} \label{unitization}
If $\fB \subseteq B(\Hil)$ is a nonunital $\bw$-algebra, then its unitization $\fB^1$ is a $\bw$-algebra in $B(\Hil)$, and for $(b_n), b \in \fB$ and $(\lambda_n),\lambda \in \C,$ we have $b_n+\lambda_n I_\Hil \xrightarrow{WOT} b + \lambda I_\Hil$ if and only if $b_n \xrightarrow{WOT} b$ and $\lambda_n \to \lambda.$
\end{prop}

\proof
If $(b_n+\lambda_n I_\Hil)$ is a sequence in $\fB^1$ converging WOT to $T$ in $B(\Hil),$ then $(\lambda_n)$ is bounded, hence has a subsequence $(\lambda_{n_k})$ converging to some $\lambda \in \C.$ So $b_{n_k} \xrightarrow{WOT} T- \lambda I_\Hil$, and thus $T \in \fB^1.$ The last claim is a short exercise using the fact that a bounded sequence in $\C$ converges iff every convergent subsequence has the same limit.
\endproof

\section{$\bw$-modules}

\begin{defn} \label{ssmod def}
A right (resp.\ left) $C^*$-module $\fX$ over a $\bw$-algebra $\fB \subseteq B(\Hil)$ is called a \emph{right (resp.\ left) $\bw$-module} if the canonical image of $\fX$ in $B(\Hil, \fX \otimes_\fB \Hil)$ (resp.\ in $B(\overline{\fX} \otimes_\fB \Hil, \Hil)$) is WOT sequentially closed. As with $C^*$-modules, ``$\fX$ is a $\bw$-module" means ``$\fX$ is a right $\bw$-module." We usually only explicitly prove results for right $\bw$-modules, but in these cases there is always an easily translated ``left version."
\end{defn}

Note the evident facts that every $\bw$-algebra is a $\bw$-module over itself (this will be generalized in Theorem \ref{pics}) and that a $\bw$-module $\fX$ over a non-unital $\bw$-algebra $\fB$ is canonically a $\bw$-module over $\fB^1$ (indeed, the algebraic module tensor products $\fX \odot_\fB \Hil$ and $\fX \odot_{\fB^1} \Hil$ coincide, so we have equality of the Hilbert spaces $\fX \otimes_\fB \Hil = \fX \otimes_{\fB^1} \Hil)$.

We will show shortly (Proposition \ref{selfdual implies bw}) that every selfdual $C^*$-module over a $\bw$-algebra is a $\bw$-module, but the converse is not true. Indeed, if $\fB$ is a nonunital $\bw$-algebra (e.g.\ the bounded Baire functions on a non-$\sigma$-compact locally compact Hausdorff space $X$, or $\Bw(\K(\Hil))$ for nonseparable $\Hil$) viewed as a $\bw$-module over itself, then $\fB$ is not selfdual since the identity map on $\fB$ is not of the form $x \mapsto y^* x$ for some $y \in \fB.$ However, we will show in Theorem \ref{bw=sd for cg} that these notions do coincide in the case of $\swb$countably generated $C^*$-modules over $\bw$-algebras.

\begin{lemma} \label{weak* conv lemma}
Let $X$ be a $C^*$-module over a $\bw$-algebra $\fB \subseteq B(\Hil)$. For a sequence $(x_n) \in X$ and $x \in X,$ we have $\langle x_n |y \rangle \xrightarrow{WOT} \langle x|y \rangle$ for all $y \in X$ if and only if $x_n \xrightarrow{WOT} x$ in $B(\Hil, X \otimes_\fB \Hil).$
\end{lemma}

\proof
($\Rightarrow$) Suppose that $\langle x_n |y \rangle \xrightarrow{WOT} \langle x|y \rangle$ for all $y \in X.$ For each $n$, let $\varphi_n: X \to \fB$ be the bounded linear map defined by $\varphi_n(y) = \langle x_n|y \rangle.$ Then for any $y \in X,$ $\sup_n \|\varphi_n(y) \| = \sup_n \| \langle x_n|y \rangle \| < \infty$ since the sequence $(\langle x_n|y \rangle)$ is WOT-convergent, hence bounded. By the uniform boundedness principle, $\sup_n \|x_n\| = \sup_n \|\varphi_n\| < \infty.$ Since $$\langle x_n(\zeta), y \otimes \eta \rangle = \langle x_n \otimes \zeta, y \otimes \eta \rangle = \langle \zeta, \langle x_n|y \rangle \eta \rangle \longrightarrow \langle \zeta, \langle x|y \rangle \eta \rangle =  \langle x_n \otimes \zeta, y \otimes \eta \rangle = \langle x(\zeta), y \otimes \eta \rangle$$ for all $\zeta, \eta \in \Hil$ and $y \in X,$ and since elements of the form $y \otimes \eta$ are total in $X \otimes_\fB \Hil,$ it follows from a triangle inequality argument that $x_n \xrightarrow{WOT} x$ in $B(\Hil, X \otimes_\fB \Hil).$

($\Leftarrow$) Assume $x_n \xrightarrow{WOT} x$ in $B(\Hil, X \otimes_\fB \Hil),$ and take $y \in X.$ Then for any $\zeta, \eta \in  \Hil,$ we have $$\langle \zeta, \langle x_n|y \rangle \eta \rangle = \langle x_n(\zeta), y \otimes \eta \rangle \to \langle x(\zeta), y \otimes \eta \rangle = \langle \zeta, \langle x|y \rangle \eta \rangle,$$ so that $\langle x_n | y \rangle \xrightarrow{WOT} \langle x|y \rangle.$
\endproof

\begin{note} \label{dist notation}
A similar result holds for left $C^*$-modules---namely, if $X$ is a left $C^*$-module over a $\bw$-algebra $\fA \subseteq B(\Kil)$ with $\fA$-valued inner product $\langle \cdot | \cdot \rangle_\fA$, then $\langle  x_n |y \rangle_\fA \xrightarrow{WOT} \langle x|y \rangle_\fA$ for all $y \in X$ if and only if $x_n \xrightarrow{WOT} x$ in $B(\overline{X} \otimes_\fA \Kil,\Kil).$ If $X$ is both a left $\bw$-module over $\fA$ and a right $\bw$-module over $\fB,$ there is thus the potential for confusion in an expression like $``x_n \xrightarrow{WOT} x.$" To distinguish, we write: $$x_n \xrightarrow{{}_\fA WOT} x \text{ iff } \langle  x_n |y \rangle_\fA \xrightarrow{WOT} \langle x|y \rangle_\fA \text{ for all } y \in X$$ and $$x_n \xrightarrow{WOT_\fB} x \text{ iff } \langle  x_n |y \rangle_\fB \xrightarrow{WOT} \langle x|y \rangle_\fB \text{ for all } y \in X$$ where $\langle \cdot | \cdot \rangle_\fA$ denotes the $\fA$-valued inner product and $\langle \cdot | \cdot \rangle_\fB$ denotes the $\fB$-valued inner product on $X$. Note that these notations make good sense even if $\fA$ and $\fB$ are concrete $C^*$-algebras that are not necessarily WOT sequentially closed.
\end{note}

The following proposition is often helpful when proving that a $C^*$-module is a $\bw$-module, and we will use it for this purpose many times.

\begin{prop} \label{bw char}
Let $\fX$ be a $C^*$-module over a $\bw$-algebra $\fB \subseteq B(\Hil)$. The following are equivalent:
\begin{enumerate}
\item[\rm{(1)}] $\fX$ is a $\bw$-module;
\item[\rm{(2)}] whenever $(x_n)$ is a sequence in $\fX$ such that $\langle x_n|y \rangle$ is WOT-convergent in $B(\Hil)$ for all $y \in \fX,$ then there is a (unique) $x \in \fX$ such that $\langle x_n|y \rangle \xrightarrow{WOT} \langle x|y \rangle$ for all $y \in \fX$;
\item[\rm{(3)}] the space $\hat{\fX}:=\{\langle x | \cdot \rangle : x \in \fX \}$ is point-WOT sequentially closed in $B_\fB(\fX,\fB).$
\end{enumerate}
\end{prop}

\proof
$(1) \implies (2).$ Assume $\fX$ is a $\bw$-module, and let $(x_n)$ be a sequence in $\fX$ such that $\langle x_n|y \rangle$ is WOT-convergent in $B(\Hil)$ for all $y \in \fX.$ Then $\langle x_n(\zeta), y \otimes \eta \rangle = \langle \zeta, \langle x_n|y \rangle \eta \rangle$ is convergent for all $\zeta, \eta \in \Hil$ and $y \in \fX.$ Since $(x_n)$ is a bounded sequence (as in the proof of the forward direction of the previous lemma), it follows that $\langle x_n(\zeta),\xi \rangle$ converges for all $\zeta \in \Hil$ and $\xi \in \fX \otimes_\fB \Hil.$ It follows from a standard argument using the correspondence betweeen operators and bounded sesquilinear maps that there is an operator $T \in B(\Hil,\fX \otimes_\fB \Hil)$ satisfying $\langle T(\zeta),\xi \rangle = \lim_n \langle x_n(\zeta),\xi \rangle$ for $\zeta \in \Hil$ and $\xi \in \fX \otimes_\fB \Hil.$ Thus $x_n \xrightarrow{WOT} T$, so by assumption, $T \in \fX.$ By the backward direction of the previous lemma, $\langle x_n|y \rangle \xrightarrow{WOT} \langle T|y \rangle$ for all $y \in \fX.$ Uniqueness follows from the usual argument that the canonical map $\fX \to B_\fB(\fX,\fB)$ is one-to-one.

$(2) \implies (1).$  Assuming (2), let $(x_n)$ be a sequence in $\fX$ such that $x_n \xrightarrow{WOT} T$ in $B(\Hil, \fX \otimes_\fB \Hil).$ Then $$\langle \zeta, \langle x_n|y \rangle \eta \rangle = \langle x_n(\zeta),y \otimes \eta \rangle \to \langle T(\zeta), y \otimes \eta \rangle$$ for all $\zeta,\eta \in \Hil$ and $y \in \fX.$ It follows that $\langle x_n|y \rangle$ is WOT-convergent for all $y \in \fX$, so by assumption there is an $x \in \fX$ such that $\langle x_n|y \rangle \xrightarrow{WOT} \langle x|y \rangle$ for all $y \in \fX.$ By the forward direction of the previous lemma, $x_n \xrightarrow{WOT} x$, so that $T=x \in \fX.$

The equivalence of (2) and (3) follows by noting that if $(x_n)$ is a sequence in $\fX$ such that $\langle x_n|y \rangle$ is WOT-convergent in $B(\Hil)$ for all $y \in \fX,$ then $y \mapsto \text{WOT-}\lim_n \langle x_n|y \rangle$ defines an operator in $B_\fB(\fX,\fB).$
\endproof

%

\begin{prop} \label{selfdual implies bw}
If $\fX$ is a selfdual $C^*$-module over a $\bw$-algebra $\fB \subseteq B(\Hil),$ then $\fX$ is a $\bw$-module.
\end{prop}

\proof
Let $(x_n)$ be a sequence in $\fX$ such that $\langle x_n|y \rangle$ is WOT-convergent for all $y \in \fX.$ Define $\psi: \fX \to \fB$ by setting $\psi(y) = \text{WOT-}\lim_n \langle x_n|y \rangle.$ It is easy to check that $\psi \in B_\fB(\fX,\fB)$, so by assumption $\psi = \langle x| \cdot \rangle$ for some $x \in \fX.$ But this means $\langle x_n |y \rangle \xrightarrow{\text{WOT}} \langle x|y \rangle$ for all $y \in \fX,$ and so by Proposition \ref{bw char}, $\fX$ is a $\bw$-module.
\endproof

One of the most basic results in the theory of $C^*$-modules (and one that is fundamental in the theory of Morita equivalence) is the fact that a right $C^*$-module $X$ over a $C^*$-algebra $A$ is a left $C^*$-module over $\K_A(X).$ Analogously, if $Y$ is a right $W^*$-module over a $W^*$-algebra $M$, then $\B_M(Y)$ is a $W^*$-module, and $Y$ is a left $W^*$-module over $\B_M(Y)$. The following proposition and theorem show that the obvious $\bw$-analogues of these statements are true. (Note that the following proposition generalizes the easy fact that the multiplier algebra and left multiplier algebra of a $\bw$-algebra are WOT sequentially closed. Indeed, in the special case $\fX=\fB$, we have by Proposition \ref{bigdiagram} that $\B_\fB(\fX) = M(\fB)$ and $B_\fB(\fB)=LM(\fB)$.)

\begin{prop} \label{mapping spaces are bw}
If $\fX$ is a right $\bw$-module over a $\bw$-algebra $\fB \subseteq B(\Hil)$, then $\mathbb{B}_\fB(\fX)$ and $B_\fB(\fX)$ are WOT sequentially closed in $B(\fX \otimes_\fB \Hil).$ For a sequence $(T_n)$ and element $T$ in $B_\fB(\fX)$, $T_n \xrightarrow{WOT} T$ if and only if $T_n(x) \xrightarrow{WOT_\fB} T(x)$ for all $x \in \fX.$
\end{prop}

\proof
Let $(T_n)$ be a sequence in $B_\fB(\fX) \subseteq B(\fX \otimes_\fB \Hil)$ such that $T_n \xrightarrow{WOT} T$ for some $T \in B(\fX \otimes_\fB \Hil).$ Then for $x,y \in \fX$ and $\zeta, \eta \in \Hil,$
$$\langle \zeta, \langle T_n(x)|y \rangle \eta \rangle = \langle T_n(x) \otimes \zeta, y \otimes \eta \rangle = \langle T_n(x \otimes \zeta), y \otimes \eta \rangle \to \langle T(x \otimes \zeta), y \otimes \eta \rangle.$$ Hence, fixing $x \in \fX,$ we have $\langle T_n(x)|y \rangle$ is WOT-convergent for all $y \in \fX$. By Proposition \ref{bw char}, there is a unique element, call it $\tilde{T}(x)$, in $\fX$ such that $$\langle T_n(x)|y \rangle \xrightarrow{WOT} \langle \tilde{T}(x)|y \rangle \text{ for all } y \in \fX.$$ Doing this for each $x \in \fX$ yields a map $\tilde{T}:\fX \to \fX.$ Since $\|\tilde{T}(x)\| = \sup\{\| \langle \tilde{T}(x)|y \rangle \|: y \in \Ball(\fX)\}$ and $\|\langle \tilde{T}(x)|y \rangle\| \leq (\sup_n \|T_n\|) \|x\| \|y\|$, we see that $\tilde{T}$ is bounded, and further direct arguments show that $\tilde{T} \in B_\fB(\fX).$ That $\tilde{T}$ coincides with $T$ in $B(\fX \otimes_\fB \Hil)$ follows by combining the two displayed expressions. Hence $B_\fB(\fX)$ is WOT sequentially closed.

Now we show that $\mathbb{B}_\fB(\fX)$ is WOT sequentially closed. If $(S_n)$ is a sequence in $\B_\fB(\fX)$ converging weakly to $S \in B(\fX \otimes_\fB \Hil),$ then by what we just proved, $S \in B_\fB(\fX)$. Since the adjoint is WOT-continuous, we also have $S^* \in B_\fB(\fX)$, where $S^*$ denotes the adjoint of $S$ as a Hilbert space operator in $B(\fX \otimes_\fB \Hil).$ For $x,y \in \fX$ and $\zeta, \eta \in \Hil,$ we have $$\langle \zeta, \langle S(x)|y \rangle \eta \rangle = \langle S(x \otimes \zeta), y \otimes \eta \rangle = \langle x \otimes \zeta, S^*(y \otimes \eta) \rangle = \langle \zeta, \langle x|S^*(y) \rangle \eta \rangle.$$ Hence $\langle S(x)|y \rangle = \langle x|S^*(y) \rangle,$ and so $S \in \B_\fB(\fX).$

For the final statement, we proved in the first paragraph above that if $T_n \xrightarrow{WOT} T$ in $B_\fB(\fX) \subseteq B(\fX \otimes_\fB \Hil)$, then $\langle T_n(x)|y \rangle \xrightarrow{WOT} \langle T(x)|y \rangle$ for all $x, y \in \fX,$ which is the same as saying $T_n(x) \xrightarrow{WOT_\fB} T(x)$ for all $x \in \fX$. Conversely, if $\langle T_n(x)|y \rangle \xrightarrow{WOT} \langle T(x)|y \rangle$ for all $x, y \in \fX,$ then $(T_n)$ is bounded by the uniform boundedness principle, and $\langle T_n(x \otimes \zeta), y \otimes \eta \rangle \to \langle T(x \otimes \zeta), y \otimes \eta \rangle$ for all $\zeta, \eta \in \Hil.$ A triangle inequality argument gives that $T_n \xrightarrow{WOT} T$ in $B(\fX \otimes_\fB \Hil).$
\endproof

Hence $\mathbb{B}_\fB(\fX)$ is a $\bw$-algebra in $B(\fX \otimes_\fB \Hil)$, and so $\Bw(\K_\fB(\fX))$, the WOT sequential closure of $\K_\fB(\fX)$ in $B(\fX \otimes_\fB \Hil)$, is contained in $\mathbb{B}_\fB(\fX)$. Since $\fX$ is a left $C^*$-module over $\B_\fB(\fX)$ with inner product taking values in $\K_\fB(\fX),$ $\fX$ is also a left $C^*$-module over the $\bw$-algebra $\Bw(\K_\fB(\fX)).$ We show in Theorem \ref{left mod over ops} that $\fX$ is in fact a $\bw$-module over $\Bw(\K_\fB(\fX)).$

We will later show (Proposition \ref{cg adj ops}), that $\Bw(\K_\fB(\fX))=\mathbb{B}_\fB(\fX)$ in the special case that $\fX$ is ``$\Sigma^*_\fB$-countably generated." We do not know of any other example (outside the $\Sigma^*_\fB$-countably generated case) in which equality holds here, but note that equality does not hold in general---for example, if $\fB$ is a nonunital $\bw$-algebra, then $\fB \cong \Bw(\K_\fB(\fB))$ is not equal to $\B_\fB(\fB)$ since the latter is unital.

\begin{lemma} \label{leftconv=rightconv}
Let $\fX$ is a right $\bw$-module over a $\bw$-algebra $\fB \subseteq B(\Hil)$. For a sequence $(x_n)$ and element $x$ in $\fX,$ 
$x_n \xrightarrow{{}_{\Bw(\K_\fB(\fX))} WOT} x$ if and only if $x_n \xrightarrow{WOT_\fB} x$.
\end{lemma}

\proof
The claim is that $|x_n \rangle \langle w| \xrightarrow{WOT} |x \rangle \langle w|$ in $B(\fX \otimes_\fB \Hil)$ for all $w \in \fX$ if and only if $\langle x_n|z \rangle \xrightarrow{WOT} \langle x|z \rangle$ in $B(\Hil)$ for all $z \in \fX.$ Assuming the former, it follows from the uniform boundedness principle that $(x_n)$ is a bounded sequence, and routine calculations give
$$\langle \langle w|y \rangle \zeta, \langle x_n|z \rangle \eta \rangle = \langle |x_n \rangle \langle w|(y \otimes \zeta),z \otimes \eta \rangle \longrightarrow \langle |x \rangle \langle w|(y \otimes \zeta),z \otimes \eta \rangle = \langle \langle w|y \rangle \zeta, \langle x|z \rangle \eta \rangle$$
for all $w,y,z \in \fX$ and $\zeta, \eta \in \Hil.$ Our usual boundedness/density arguments show that if $P \in B(\Hil)$ is the projection onto the closed subspace of $\Hil$ generated by $\{ \langle x|y \rangle \zeta : x,y \in \fX \text{ and } \zeta \in \Hil \}$, then for any $\xi, \eta \in \Hil$ and $z \in \fX,$ we have
$$\langle \xi, \langle x_n|z \rangle \eta \rangle = \langle P \xi, \langle x_n|z \rangle \eta \rangle \longrightarrow \langle P \xi, \langle x|z \rangle \eta \rangle = \langle \xi, \langle x|z \rangle \eta \rangle.$$ Hence $\langle x_n|z \rangle \xrightarrow{WOT} \langle x|z \rangle$ in $B(\Hil)$ for all $z \in \fX.$  The converse is similar.
\endproof

\begin{thm} \label{left mod over ops}
If $\fX$ is a right $\bw$-module over a $\bw$-algebra $\fB \subseteq B(\Hil)$, then $\fX$ is a left $\bw$-module over the $\bw$-algebra $\Bw(\K_\fB(\fX)) \subseteq B(\fX \otimes_\fB \Hil).$
\end{thm}

\proof
By the ``left version" of Proposition \ref{bw char}, we need to show that if $(x_n)$ is a sequence in $\fX$ such that $|x_n \rangle \langle y|$ is WOT-convergent in $B(\fX \otimes_\fB \Hil)$ for all $y \in \fX,$ then there is an $x \in \fX$ such that $|x_n \rangle \langle y| \xrightarrow{WOT} |x \rangle \langle y|$ for all $y \in \fX.$ If $|x_n \rangle \langle y|$ is WOT-convergent in $B(\fX \otimes_\fB \Hil)$ for all $y \in \fX,$ then arguments from the first paragraph of the proof of Lemma \ref{leftconv=rightconv} show that $\langle x_n|z \rangle$ is WOT-convergent for all $z \in \fX.$ By Proposition \ref{bw char}, there is an $x \in \fX$ such that $\langle x_n|z \rangle \to \langle x|z \rangle$ for all $z \in \fX,$ and by Lemma \ref{leftconv=rightconv}, $|x_n \rangle \langle y| \xrightarrow{WOT} |x \rangle \langle y|$ for all $y \in \fX.$
\endproof

A \emph{ternary ring of operators} (abbreviated \emph{TRO}) is a subspace $Z \subseteq B(\Hil,\Kil)$, for Hilbert spaces $\Hil,\Kil$, such that $x y^* z \in Z$ for all $x,y,z \in Z$; and a \emph{corner} of a $C^*$-algebra $A$ is a subspace of the form $pAq$ for projections $p,q \in M(A)$. (This is slightly different from the usual definition of a corner as a subspace of the form $pAp^\perp,$ but every corner in our sense can be identified with a corner in the usual sense of a different $C^*$-algebra, so the two definitions are not essentially different.) Note that if $Z$ is a TRO in $B(\Hil,\Kil)$, then there is a canonical triple isomorphism (see \cite[8.3.1]{BLM}) identifying $Z$ with a TRO in $B(\Hil,[Z \Hil])$. So, just as for $C^*$-algebras, there is no real loss in assuming from the outset that a TRO is \emph{nondegenerate}, i.e.\ that $[Z \Hil]=\Kil.$

In analogy with the situation in $C^*$-module theory and $W^*$-module theory, $\bw$-modules are essentially the same as WOT sequentially closed TROs, and essentially the same as corners of $\bw$-algebras. The next theorem gives the details for how to move from one of these ``pictures" to another. To prepare for this, we first describe the $\bw$-version of the ``linking algebra" of a $C^*$-module.

\begin{prop} \label{bw linking}
If $\fX$ is a $\bw$-module over a $\bw$-algebra $\fB \subseteq B(\Hil)$, then $\mathcal{L}^\Bw(\fX):=\left[ \begin{matrix} \Bw(\K_\fB(\fX)) & \fX \\ \overline{\fX} & \fB \end{matrix} \right]$ is a $\bw$-algebra in $B((\fX \otimes_\fB \Hil) \oplus^2 \Hil)$.
\end{prop}

\proof
It is very easy to show that a sequence of $2 \times 2$ matrices in $\mathcal{L}^\Bw(\fX)$ converges WOT to a $2 \times 2$ matrix $\xi \in B((\fX \otimes_\fB \Hil) \oplus^2 \Hil)$ if and only if each of the entries converges WOT to the corresponding entry in $\xi.$ Since each of the four corners of $\mathcal{L}^\Bw(\fX)$ is WOT sequentially closed, the result follows.
\endproof


In the following theorem, when we say $``\fX \cong (1-p)\fC p$ and $\fB \cong p \fC p$ under isomorphisms preserving all the $\bw$-module structure," we mean that there is an isometric isomorphism $\varphi: \fX \to (1-p)\fC p$ and a $\bw$-isomorphism $\psi: \fB \to p \fC p$ such that $\varphi(xb)=\varphi(x)\psi(b)$ and $\psi(\langle x|y \rangle) = \varphi(x)^*\varphi(y)$ for all $x,y \in \fX$ and $b \in \fB.$ Note that these conditions imply that $\varphi(x_n) \xrightarrow{WOT_{p \fC p}} \varphi(x)$ whenever $x_n \xrightarrow{WOT_\fB} x.$

\begin{thm} \label{pics}

\begin{enumerate}
\item[\rm{(1)}] If $\fX$ is a $\bw$-module over $\fB \subseteq B(\Hil)$, then $\fX$ is a WOT sequentially closed TRO in $B(\Hil, \fX \otimes \Hil).$ Conversely, if $\fZ$ is a nondegenerate WOT sequentially closed TRO in $B(\Kil_1,\Kil_2),$ then $\fZ$ is a $\bw$-module over $\Bw(\fZ^*\fZ)$ with the obvious module action and inner product $\langle z_1|z_2 \rangle = z_1^* z_2.$

\item[\rm{(2)}]  If $\fY = p\fD q$ is a corner of a $\bw$-algebra $\fD,$ then $\fY$ is canonically a $\bw$-module over $q \fD q$. Conversely, if $\fX$ is a $\bw$-module over $\fB \subseteq B(\Hil)$, then there exists a $\bw$-algebra $\fC \subseteq B(\Kil)$ and a projection $p \in M(\fC)$ such that $\fX \cong (1-p)\fC p$ and $\fB \cong p\fC p$ under isomorphisms preserving all the $\bw$-module structure.
\end{enumerate}

\end{thm}

\proof
(1) The forward direction follows immediately from the definition of $\bw$-modules. For the converse, we must first show that $\fZ$ is closed under right multiplication by elements in $\Bw(\fZ^* \fZ).$ Fixing $z \in \fZ,$ the set $\mathscr{S}_z=\{b \in \Bw(\fZ^* \fZ) : zb \in \fZ\}$ contains $\fZ^* \fZ$ since $\fZ$ is a TRO, and an easy argument shows that $\mathscr{S}_z$ is WOT sequentially closed, so that $\mathscr{S}_z = \Bw(\fZ^* \fZ).$ So $\fZ$ is a right module over $\Bw(\fZ^* \fZ)$, and it is straightforward to show that it is a $C^*$-module over $\Bw(\fZ^* \fZ)$ with the canonical inner product. To prove that $\fZ$ is a $\bw$-module, note that under the canonical unitary $\fZ \otimes_{\Bw(\fZ^* \fZ)} \Kil_1 \cong [\fZ \Kil_1]=\Kil_2$, the embedding $\fZ \hookrightarrow B(\Kil_1,\fZ \otimes_{\Bw(\fZ^* \fZ)} \Kil_1)$ coincides with the inclusion $\fZ \subseteq B(\Kil_1,\Kil_2),$ so it follows from the definition that $\fZ$ is a $\bw$-module over $\Bw(\fZ^* \fZ).$


(2) For the forward direction, first note the easy fact that if $\fD$ is a $\bw$-algebra in $B(\Kil)$ and $q \in M(\fD) \subseteq B(\Kil)$ is a projection, then $q\fD q$ is a $\bw$-algebra in $B(q \Kil).$ Showing that $\fY$ is a $\bw$-module over $q \fD q$ is then a short exercise either using the definition as in the proof of (1) or employing Proposition \ref{bw char}. The converse follows from Proposition \ref{bw linking} with $\fC= \mathcal{L}^\Bw(\fX)$ and $p= \left[ \begin{matrix} 0 & 0 \\ 0 & 1 \end{matrix} \right].$
\endproof


It is an interesting and useful fact that a $C^*$-module over a $W^*$-algebra always admits a ``selfdual completion," that is, a unique $W^*$-module containing the original $C^*$-module as a weak*-dense submodule. Hamana in \cite{Hamana} and Lin in \cite{Lin92} also proved that a $C^*$-module $X$ over a monotone complete $C^*$-algebra admits a selfdual completion, and Hamana proved uniqueness under the condition that $X^\perp = (0)$. The proposition below gives existence of a ``$\bw$-module completion" analogous to the selfdual completion.

Note that an easy modification of Lemma \ref{completion lem} and Proposition \ref{completion} gives another proof of the existence of the selfdual completion of a $C^*$-module over a $W^*$-algebra (this is surely known to experts though).

For a $C^*$-module $X$ over a nondegenerate $C^*$-algebra $\fB \subseteq B(\Hil),$ recall a few canonical embeddings from Proposition \ref{bigdiagram}: 
$$X \cong \K_\fB(\fB,X) \hookrightarrow B(\Hil,X \otimes_\fB \Hil)$$
$$B_\fB(X,\fB) \hookrightarrow B(X \otimes_\fB \Hil, \Hil).$$
In the following lemma, the definition of $\mathscr{S}$ implicitly uses the latter, and the last few statements use the former.

\begin{lemma} \label{completion lem}
If $X$ is a $C^*$-module over a $\bw$-algebra $\fB \subseteq B(\Hil),$ then $$\mathscr{S} := \{T \in B(\Hil, X \otimes_\fB \Hil) : T^* \in B_\fB(X, \fB)\}$$ is WOT sequentially closed in $B(\Hil,X \otimes_\fB \Hil)$ and contains $X$. Hence we may view $$X \subseteq \Bw(X) \subseteq \mathscr{S} \subseteq B(\Hil,X \otimes_\fB \Hil),$$ where by $\Bw(X)$ we mean the WOT sequential closure of $X$ in $B(\Hil,X \otimes_\fB \Hil)$.
\end{lemma}

\proof
Suppose that $(T_n)$ is a sequence in $\mathscr{S}$ and $T \in B(\Hil, X \otimes_\fB \Hil)$ with $T_n \xrightarrow{WOT} T$ in $B(\Hil, X \otimes_\fB \Hil)$. Then $T_n^* \xrightarrow{WOT} T^*$ in $B(X \otimes_\fB \Hil,\Hil),$ so $$\langle T_n^*(x)\zeta,\eta \rangle = \langle T_n^*(x \otimes \zeta), \eta \rangle \longrightarrow \langle T^*(x \otimes \zeta),\eta \rangle \text{ for all } x \in X \text{ and } \zeta,\eta \in \Hil,$$ and it follows that for each $x \in X,$ the sequence $(T_n^*(x))$ converges WOT in $\fB \subseteq B(\Hil).$ Define a map $\tau: X \to \fB$ by $\tau(x) = \text{WOT-} \lim_n T_n^*(x).$ It is direct to check that $\tau$ is in $B_\fB(X,\fB),$ and since $\langle T^*(x \otimes \zeta),\eta \rangle = \langle \tau(x) \zeta, \eta \rangle = \langle \tau(x \otimes \zeta),\eta \rangle$ for all $x \in X$ and $\zeta, \eta \in \Hil,$ we may conclude that $T^*$ coincides with $\tau$ under the embedding $B_\fB(X,\fB) \hookrightarrow B(X \otimes_\fB \Hil, \Hil).$ Hence $T \in \mathscr{S}$, and so $\mathscr{S}$ is WOT sequentially closed. All the other claims are evident.
\endproof

Note that it follows quickly from the definitions that $\xp=X$ if and only if $X$ is a $\bw$-module. To see that $\xp$ and $\mathscr{S}$ may be different, take $X = \fB$ for a nonunital $\bw$-algebra $\fB$. Then $\xp = \fB \not= \mathscr{S}$ since $I_\Hil \in \mathscr{S}$.

To explain some terminology that appears in the following theorem and later on in this paper, a $C^*$-submodule $X$ of a $\bw$-module $\fX$ is said to be \emph{WOT$_\fB$ sequentially dense} if $\fX$ is the only subset of itself that contains $X$ and is closed under limits of WOT$_\fB$-convergent sequences. Note that this may be different from saying that every element in $\fX$ is the WOT$_\fB$-limit of a sequence in $X.$

\begin{thm} \label{completion}
If $X$ is a $C^*$-module over a $\bw$-algebra $\fB \subseteq B(\Hil),$ then in the notation of the preceding lemma, $\Bw(X)$ has a $\fB$-valued inner product making it into a $\bw$-module that contains $X$ as a WOT$_\fB$ sequentially dense submodule and has $$\langle \tau|x \rangle = \tau^*(x) \text{ for all } \tau \in \xp, x \in X.$$ Moreover, the operator norm $\xp$ inherits from $B(\Hil,X \otimes_\fB \Hil)$ coincides with this $C^*$-module norm.
\end{thm}

\proof
We first show that $\xp$ is a right $\fB$-module with the canonical module action coming from the inclusions $\xp \subseteq B(\Hil,X \otimes_\fB \Hil)$ and $\fB \subseteq B(\Hil).$ Fix $b \in \fB,$ and let $\mathscr{S}=\{x \in \xp : xb \in \xp\}.$ Then $\mathscr{S}$ is WOT sequentially closed and contains $X$, so $\mathscr{S}=\xp.$ Since $b \in \fB$ was arbitrary, we have shown that $xb \in \xp$ for all $x \in \xp$ and $b \in \fB.$

Now note that for any $K,L \in \mathbb{K}_\fB(\fB,X)$, $K^*L$ is in $\mathbb{K}_\fB(\fB)=\fB \subseteq B(\Hil).$ Using this, arguments of the sort used in the previous paragraph (or in Proposition \ref{closure is $C^*$}) shows that $S^*T \in \fB$ for all $S, T \in \xp$. Define a $\fB$-valued inner product on $\xp$ by $\langle S|T \rangle := S^*T.$ With this inner product and the right $\fB$-module structure it inherits from $B(\Hil,X \otimes_\fB \Hil)$, it is easy to check that $\xp$ is a $C^*$-module over $\fB.$ It is also straightforward to check that the centered equation in the claim holds.

To see that $\xp$ is a $\bw$-module, suppose that $\tau_n$ is a sequence in $\xp$ such that $\langle \tau_n|\sigma \rangle$ converges WOT in $\fB$ for all $\sigma \in \xp.$ In particular, $\langle \tau_n|x \rangle = \tau_n^*(x)$ converges WOT to an element in $\fB$, call it $\tau^*(x)$, for each $x \in X.$ Routine arguments show that $\tau^*: X \to \fB$ thus defined is in $B_\fB(X,\fB)$ and that $\tau_n^* \xrightarrow{WOT} \tau^*$ in $B(X \otimes_\fB \Hil,\Hil),$ so that $\langle \tau_n|\sigma \rangle = \tau_n^* \sigma \xrightarrow{WOT} \tau^* \sigma = \langle \tau|\sigma \rangle$ in $\fB$ for all $\sigma \in \xp.$ Hence $\xp$ is a $\bw$-module by Proposition \ref{bw char}.

Note that we have demonstrated that WOT$_\fB$-convergence of a sequence in $\xp$ is the same as WOT-convergence in $\xp$ considered as a subset of $B(\Hil, X \otimes_\fB \Hil).$ This fact combined with the definition of $\xp$ gives that $X$ is WOT$_\fB$ sequentially dense in $\xp.$

The last claim follows immediately from the definition of the inner product: $\| \langle \tau|\tau \rangle \|_{\xp}^2 = \| \tau^* \tau\|_{B(\Hil)} = \|\tau\|_{B(\Hil,X \otimes_\fB \Hil)}^2.$
\endproof

Unfortunately, we were not able in general to prove uniqueness of the above construction with conditions as simple as those in the $W^*$-case or monotone complete case (but see Proposition \ref{unique compl} for a special case).

\begin{defn}For a $C^*$-module $X$ over a $\bw$-algebra $\fB,$ a \emph{$\bw$-module completion} of $X$ is any $\bw$-module $\fX$ over $\fB \subseteq B(\Hil)$ such that:
\begin{enumerate}
\item[\rm{(1)}]  $\fX$ contains $X$ as a WOT$_\fB$ sequentially dense submodule;
\item[\rm{(2)}]  the $\fB$-valued inner product on $\fX$ extends that of $X$;
\item[\rm{(3)}]  $\|\xi\| = \sup\{\| \langle \xi|x \rangle \| : x \in \text{Ball}(X)\}$ for all $\xi \in \fX$;
\item[\rm{(4)}]  if $(\xi_n)$ is a sequence in $\fX$ such that $(\langle \xi_n|x \rangle)$ is WOT-convergent for all $x \in X,$ then there is a $\xi \in \fX$ such that $\xi_n \xrightarrow{WOT_\fB} \xi.$
\end{enumerate}
\end{defn}

\begin{prop} \label{uniq}
If $X$ is a $C^*$-module over a $\bw$-algebra $\fB \subseteq B(\Hil),$ then the $\bw$-module $\xp$ of the previous theorem is the unique $\bw$-module completion of $X$ (up to unitary isomorphism).
\end{prop}

\proof
It follows immediately from the previous theorem that $\xp$ satisfies (1) and (2) in the definition of a $\bw$-module completion. To see (3), for $\tau \in \xp,$ we have
\begin{equation*}
\begin{split}
\|\tau\|_{\xp} &= \|\tau\|_{B(\Hil,X \otimes_\fB \Hil)} \\ &= \|\tau^*\|_{B(X \otimes_\fB \Hil,\Hil)}\\ &= \|\tau^*\|_{B_\fB(X,\fB)} \\
&=\sup\{ \|\tau^*(x)\| : x \in \Ball(X)\} \\
&= \sup\{ \| \langle \tau|x \rangle \|: x \in \Ball(X)\}.
\end{split}
\end{equation*}
(The first equality here follows from the last claim in Theorem \ref{completion}, the third equality follows from Lemma \ref{completion lem} and the isometric embedding $B_\fB(X,\fB) \hookrightarrow B(X \otimes_\fB \Hil,\Hil),$ and the final equality follows from the centered equation in Theorem \ref{completion}.) The argument for (4) basically follows the second paragraph of the proof of the previous theorem.

To prove uniqueness, suppose that $\fY$ is another $\bw$-module completion of $X$, and denote its $\fB$-valued inner product by $(\cdot|\cdot)$. Define maps $V: \fY \to B_\fB(X,\fB)$ and $U: \fY \to B(\Hil, X \otimes_\fB \Hil)$ by $V(\xi)(x)= ( \xi|x )$ and $U(\xi)=V(\xi)^*$ for $\xi \in \fY$ and $x \in X.$ We will show that $U$ is a $\fB$-linear isometry with range equal to $\xp,$ and so the result follows Lance's result that every isometric, surjective module map between $C^*$-modules is a unitary (\cite[Theorem 3.5]{Lan}).

Note first the formula $$\langle U(\xi) \eta, x \otimes \zeta \rangle = \langle \eta, V(\xi)(x \otimes \zeta) \rangle = \langle \eta, V(\xi)(x) \zeta \rangle= \langle \eta, (\xi|x)\zeta \rangle$$ for $\xi \in \fY, x \in X,$ and $\eta,\zeta \in \Hil.$ An easy calculation from this shows that $U$ is linear and $\fB$-linear.

By this, if $z \in X \subseteq \fY,$ then $\langle U(z) \eta, x \otimes \zeta \rangle = \langle \eta, (z|x) \zeta \rangle = \langle z \otimes \eta, x \otimes \zeta \rangle$ for all $\zeta, \eta \in \Hil$ and $x \in X$. Hence $U(z)=z$ in $B(\Hil, X \otimes_\fB \Hil)$, and we have shown that $X \subseteq U(\fY).$

Let $\mathscr{T}=\{\xi \in \fY : U(\xi) \in \xp\}.$ We just showed that $X \subseteq \mathscr{T}$, so if we can show that $\mathscr{T}$ is WOT$_\fB$ sequentially closed, it will follow by sequential WOT$_\fB$-density of $X$ in $\fY$ that $\mathscr{T}=\fY.$ To this end, suppose that $(\xi_n)$ is a sequence in $\mathscr{T}$ with $\xi_n \xrightarrow{WOT_\fB} \xi$ in $\fY.$ By the centered line above, $$\langle U(\xi_n) \eta, x \otimes \zeta \rangle = \langle \eta, (\xi_n|x) \zeta \rangle \longrightarrow \langle \eta, (\xi|x) \eta \rangle = \langle U(\xi)\eta, x \otimes \zeta \rangle$$ for all $x \in X$ and $\eta, \zeta \in \Hil,$ so that $U(\xi_n) \xrightarrow{WOT} U(\xi)$ in $B(\Hil, X \otimes_\fB \Hil).$ Thus $U(\xi) \in \xp,$ and $\xi \in \mathscr{T}.$ So we may conclude that $\mathscr{T}=\fY,$ which is to say $U(\fY) \subseteq \xp.$

Combining the previous two paragraphs, we have $X \subseteq U(\fY) \subseteq \xp.$ So to have $U(\fY)=\xp,$ it remains to prove that $U(\fY)$ is WOT$_\fB$ sequentially closed in $\xp.$ Suppose that $(\xi_n)$ is a sequence in $\fY$ with $U(\xi_n) \xrightarrow{WOT_\fB} \tau$ in $\xp.$ Then for $\zeta, \eta \in \Hil$ and $x \in X,$ $$\langle \zeta, ( \xi_n |x ) \eta \rangle = \langle U(\xi_n)(\zeta),x \otimes \eta \rangle \to \langle \tau(\zeta),x \otimes \eta \rangle = \langle \zeta, \tau^*(x)(\eta) \rangle = \langle \zeta, \langle \tau|x \rangle \eta \rangle,$$ so that $(\xi_n|x) \xrightarrow{WOT} \langle \tau|x \rangle$ in $\fB.$ By assumption (4) in the definition above the proposition, $\xi_n \xrightarrow{WOT_\fB} \xi$ for some $\xi \in \fY,$ and the argument in the previous paragraph shows that $U(\xi_n) \xrightarrow{WOT_\fB} U(\xi),$ so that $\tau = U(\xi) \in U(\fY).$

It remains to show that $U$ is isometric. If $\xi \in \fY,$ $x \in X,$ and $\eta,\zeta \in \Hil,$ then $$\langle \langle U(\xi)|x \rangle \zeta, \eta \rangle = \langle x \otimes \zeta, U(\xi)(\eta) \rangle = \langle V(\xi)(x \otimes \zeta),\eta \rangle = \langle (\xi|x) \zeta,\eta \rangle,$$ which gives that $\langle U(\xi)|x \rangle = (\xi|x)$ for all $x \in X.$ That $\|U(\xi)\|=\|\xi\|$ now follows from assumption (3) in the definition above.
\endproof

To close this section, we provide a result that is used in the next section and seems interesting when one dwells upon the similarities between $\bw$-modules and $W^*$-modules. For $W^*$-modules $Y$ and $Z$ over $M$, we have $\{\langle y | \cdot \rangle : y \in Y\} = B_M(Y,M)$ and $\mathbb{B}_M(Y,Z)=B_M(Y,Z)$ and all the maps in both of these spaces are weak*-continuous. The following result is a $\bw$-analogue of this fact, but with an additional condition that may be taken as a weak type of the assumption of being ``countably generated" (indeed, we will see in the next section that all $\swb$countably generated $\bw$-modules meet this condition). This condition cannot be removed in general---for example, if $\fB$ is a nonunital $\bw$-algebra considered as a $\bw$-module over itself, then $\text{id}_\fB$ is in the latter set in (1) below, but is not in the former.

\begin{prop} \label{char of adj ops}
If $\fX$ is a $\bw$-module over a $\bw$-algebra $\fB \subseteq B(\Hil)$ such that $\Bw(\K_\fB(\fX))=\B_\fB(\fX),$ and $\fY$ is any other $\bw$-module over $\fB,$ then
\begin{enumerate}
\item[\rm{(1)}]  $\{\langle z| \cdot \rangle : z \in \fX\} = \{\xi \in B_\fB(\fX,\fB) : x_n \xrightarrow{WOT_\fB} x \implies \xi(x_n) \xrightarrow{WOT} \xi(x)\};$
\item[\rm{(2)}]   $\B_\fB(\fX,\fY) = \{T \in B_\fB(\fX,\fY) :  x_n \xrightarrow{WOT_\fB} x \implies T(x_n) \xrightarrow{WOT_\fB} T(x)\}.$
\end{enumerate}
\end{prop}

\proof
The forward inclusion of (1) is evident from the definitions. For the other inclusion, fix a $\xi$ as in the latter set in (1). Note that the condition $\Bw(\K_\fB(\fX))=\B_\fB(\fX)$ is equivalent to saying that $\B_\fB(\fX)$ is generated as a $\bw$-algebra by the ``finite-rank operators," that is, operators of the form $\sum_{i=1}^n|x_i \rangle \langle y_i|$ for $x_i,y_i \in \fX.$ Let $$\mathscr{T} = \{T \in \B_\fB(\fX): \xi \circ T = \langle z| \cdot \rangle \text{ and } \xi \circ T^* = \langle w| \cdot \rangle \text{ for some } z,w \in \fX\}.$$ To see that $\mathscr{T}$ is WOT sequentially closed, suppose that $(T_n)$ is a sequence in $\mathscr{T}$ with $T_n \xrightarrow{WOT} T$ in $\B_\fB(\fX) \subseteq B(\fX \otimes_\fB \Hil).$ For $x \in \fX,$ by Proposition \ref{mapping spaces are bw} we have $T_n(x) \xrightarrow{WOT_\fB} T(x)$, so that $\xi \circ T_n(x) \xrightarrow{WOT} \xi \circ T(x)$ by the assumption on $\xi$. Writing $\xi \circ T_n(x) = \langle z_n|x \rangle,$ we may conclude by  Proposition \ref{bw char} that $\xi \circ T = \langle z| \cdot \rangle$ for some $z \in \fX.$ Since $T_n \xrightarrow{WOT} T$ implies $T_n^* \xrightarrow{WOT} T^*,$ the same argument shows that $\xi \circ T^* = \langle w| \cdot \rangle$ for some $w \in \fX.$ So $\mathscr{T}$ is WOT sequentially closed. It is easy to check that $\mathscr{T}$ is a $*$-subalgebra of $\B_\fB(\fX)$ containing all the finite-rank operators; hence $\mathscr{T} = \B_\fB(\fX).$ Since $I \in \B_\fB(\fX),$ we conclude that $\xi = \langle z| \cdot \rangle$ for some $z \in \fX.$


For the forward inclusion of (2), let $S \in \B_\fB(\fX,\fY),$ and suppose $x_n \xrightarrow{WOT_\fB} x$ in $\fX.$ Then $$\langle S(x_n)|y \rangle = \langle x_n|S^*(y) \rangle \xrightarrow{WOT} \langle x|S^*(y) \rangle = \langle S(x)|y \rangle$$ for all $y \in \fX,$ so that $S(x_n) \xrightarrow{WOT_\fB} S(x).$ For the other inclusion, suppose that $T$ is in the latter set in (2). Then for any $y \in \fY,$ the map $\langle y|T(\cdot) \rangle$ is in the latter set in (1), so there is a $z \in \fX$ such that $\langle z|x \rangle = \langle y|T(x) \rangle$ for all $x \in \fX.$ Hence $T$ is adjointable.
\endproof

\begin{note}
In principle, one could work out analogous theories to that presented above for many different classes of $C^*$-algebras. For example, one could define a \emph{Borel module} to be a $C^*$-module $\fX$ over a Borel $*$-algebra $\fB \subseteq B(\Hil)$ such that $\left[ \begin{matrix} \Bw^m(\K_\fB(\fX)) & \fX \\ \overline{\fX} & \fB \end{matrix} \right]$ is monotone sequentially closed in $B((\fX \otimes_\fB \Hil) \oplus^2 \Hil),$ where $\Bw^m(\cdot)$ denotes monotone sequential closure. It would be interesting to try to work out the appropriate Borel analogues of the results for $\bw$-modules we have proved, but it does not seem clear how to do this even for the first few of our results.
\end{note}

\section{Countably generated $\bw$-modules}
Many of the most interesting results in $C^*$-module theory require some type of ``smallness" condition on either the module or the coefficient $C^*$-algebra, e.g.\ that the module is countably generated or that the $C^*$-algebra is separable or $\sigma$-unital. In this section, we study a weak sequential analogue of the condition of being a countably generated module. The more elegant results we obtain in this section indicate that $\bw$-modules meeting this countably generated condition are more similar to $W^*$-modules than are general $\bw$-modules.

The main highlights of this section are Proposition \ref{wcg} (which is analogous to some well-known equivalent conditions to being a (norm) countably generated $C^*$-module (see \cite[8.2.5]{BLM})), Theorem \ref{bw=sd for cg} (which says that in the class of ``$\bw_\fB$-sequentially countably generated" $C^*$-modules over $\bw$-algebras, the $\bw$-modules coincide with the selfdual $C^*$-modules); Proposition \ref{dpb} (an interesting result about column spaces donated to us by David Blecher); Theorem \ref{stabilization} (our analogue of Kasparov's stabilization theorem); and Proposition \ref{orth compl}.

\begin{defn}
\begin{itemize}
\item
A right $C^*$-module $X$ over a $C^*$-algebra $A \subseteq B(\Hil)$ is \emph{$\swa$countably generated} if there is a countable set $\{x_i\}_{i=1}^\infty$ such that $X$ is generated as a relatively WOT$_A$ sequentially closed subset of itself by $\{\sum_{i=1}^N x_i b_i : b_i \in A, N \in \N\}$.
\item
For a nondegenerate $C^*$-algebra $\mathfrak{C} \subseteq B(\mathcal{K}),$ say that a sequence $(e_n)$ in $\Ball(\mathfrak{C})$ is a \emph{sequential weak cai for $\mathfrak{C}$} if $e_n c \xrightarrow{\text{WOT}} c$ for all $c \in \mathfrak{C}$. Since $(e_n)$ is bounded and $[\mathfrak{C}\Kil]=\Kil,$ a triangle inequality argument shows that $(e_n)$ is a sequential weak cai if and only if $e_n \xrightarrow{\text{WOT}} I_\Kil.$ So in this case we also have $c e_n \xrightarrow{\text{WOT}} c$ for all $c \in \mathfrak{C}.$

\item A $\bw$-algebra $\fB \subseteq B(\Hil)$ is \emph{$\bw$-countably generated} (resp.\ \emph{$\bw$-singly generated}) if there is a countable (resp.\ singleton) subset $B$ of $\fB$ such that $B$ generates $\fB$ as a $\bw$-algebra, that is, $\Bw(C^*(B))=\fB$, where $C^*(B)$ is the $C^*$-algebra generated by $B$ and $\Bw(\cdot)$ denotes WOT sequential closure (see notation above Proposition \ref{closure is $C^*$}).
\end{itemize}
\end{defn}

\begin{example} \label{cg examples}
\begin{enumerate}
\item[\rm{(1)}]  If $\fC \subseteq B(\Kil)$ is a $\bw$-countably generated $\bw$-algebra (e.g.\ the $\bw$-envelope of a separable $C^*$-algebra), and $p \in M(\fC),$ then $(1-p) \fC p$ is a right $\bw$-module over $p \fC p$ (see Theorem \ref{pics} (2)), and $(1-p) \fC p$ is $\Sigma^*_{p \fC p}$-countably generated. Indeed, one may deduce this quickly from the following observation (which uses and is analogous to the fact that countably generated $C^*$-algebras are separable): if a $\bw$-algebra $\fB$ is $\bw$-countably generated, then there is a countable subset $D$ of $\fB$ such that $\Bw(D)=\fB$.

\item[\rm{(2)}]  It is immediate that if a $\bw$-algebra $\fB$, considered as a $\bw$-module over itself, is unital, then it is $\swb$countably generated. We will show in Corollary \ref{cg is unital} that the converse of this is also true.

\item[\rm{(3)}]  For a unital $\bw$-algebra $\fB,$ the column $\bw$-module $C^w(\fB)$ described above Corollary \ref{col sd} is $\swb$countably generated.
\end{enumerate}
\end{example}

The von Neumann algebra analogue of the following proposition is well-known, and since the spectral theorem still holds in $\bw$-algebras (by Proposition \ref{sp proj}), the proof is virtually the same. We thank David Blecher for pointing this result out, and for the example following.

\begin{prop} \label{comm cg is sg}
If $\fB$ is a $\bw$-countably generated commutative $\bw$-algebra, then it is $\bw$-singly generated by a selfadjoint element.
\end{prop}

\begin{note}
Related to (2) in Example \ref{cg examples}, it is easy to see that every $\bw$-countably generated $\bw$-algebra is unital (since a countable subset will generate a $\sigma$-unital, WOT sequentially dense $C^*$-subalgebra), but the converse is not necessarily true. Take for example the von Neumann algebra $\ell^\infty(I) \subseteq B(\ell^2(I))$ for a set $I$ with cardinality strictly greater than that of $\R.$ If $\ell^\infty(I)$ were $\bw$-countably generated, then by Proposition \ref{comm cg is sg} it would be $\bw$-singly generated by a selfadjoint element $x=(x_i)_{i \in I}.$ However, as the map $I \to \R,$ $i \mapsto x_i,$ cannot be one-to-one, there must be $k,l \in I,$ $k \not=l,$ with $x_k=x_l.$ Since the set $\mathscr{S}$ of $(y_i)_{i \in I}$ in $\ell^\infty(I)$ such that $y_k=y_l$ is WOT sequentially closed in $B(\ell^2(I))$ and contains $x$, we have the contradiction $\ell^\infty(I) \subseteq \mathscr{S}.$
\end{note}


The following simple lemma is a weak sequential version of some well-known characterizations of $\sigma$-unital $C^*$-algebras (cf.\ \cite[3.10.5]{Ped}).

\begin{lemma} \label{wsp lem}
If $A \subseteq B(\Hil)$ is a $C^*$-algebra, the following are equivalent:
\begin{enumerate}
\item[\rm{(1)}]  $A$ has an element $a$ such that $\psi(a)>0$ for all nonzero WOT sequentially continuous positive functionals $\psi$ on $\Bw(A)$;

\item[\rm{(2)}]  $A$ has a positive element $a$ such that $\overline{a(\Hil)} = \Hil$;

\item[\rm{(3)}]  $A$ has a positive increasing sequential weak cai.
\end{enumerate}
\end{lemma}

\proof
$(1) \implies (2).$ Let $a \in A$ be as in (1). Then $\langle a \zeta, \zeta \rangle >0$ for all nonzero $\zeta \in \Hil$, and hence $\text{Ker}(a) = \text{Ran}(a)^\perp = (0) \implies \overline{\text{Ran}(a)} = \Hil$.

$(2) \implies (3).$ Assume (2), and set $e_n = a^{1/n}.$ Then $e_n \nearrow s(a)$ in $B(\Hil),$ where $s(a)$ denotes the support projection of $a$. Since $a$ has dense range, $s(a)$ is the identity operator on $\Hil.$ So $e_n \xrightarrow{\text{WOT}} I_\Hil.$

$(3) \implies (1).$ Let $(e_n)$ be a positive increasing sequential weak cai in $A$. As mentioned in the definition, this means that $e_n \nearrow I_\Hil.$ Set $a=\sum_{n=1}^\infty 2^{-n} e_n.$ If $\psi$ is a WOT sequentially continuous positive functional on $\Bw(A)$ with $\psi(a)=0,$ then $\psi(e_n)=0$ for all $n$ since $e_n \leq a.$ But since $\psi(e_n) \nearrow \psi(I_\Hil) = \|\psi\|,$ we have $\psi=0.$
\endproof

\begin{lemma} \label{gen set}
If $X$ is a right $C^*$-module over a $\bw$-algebra $\fB \subseteq B(\Hil)$ that is $\swb$countably generated by a subset $\{x_i\},$ then $\Bw(\{\sum_{i,j=1}^n |x_i b_{ij} \rangle \langle x_j| : b_{ij} \in \fB\}) = \Bw(\mathbb{K}_\fB(X)).$
\end{lemma}

\proof
Set $\mathcal{A} := \{\sum_{i,j=1}^n |x_i b_{ij} \rangle \langle x_j| : b_{ij} \in \fB, n \in \N\}.$ Clearly $\Bw(\mathcal{A}) \subseteq \Bw(\K_\fB(X)).$ It is easy to check that $\mathcal{A}$ is a $*$-subalgebra of $\B_\fB(X)$, so $\Bw(\mathcal{A}) = \Bw(\overline{\mathcal{A}})$ is a $C^*$-algebra, and thus the inclusion $\Bw(\K_\fB(X)) \subseteq \Bw(\mathcal{A})$ will follow if we can show that $|x \rangle \langle z| \in \Bw(\mathcal{A})$ for all $x, z \in X.$ Fix $k \in \mathbb{N}$ and $b \in \fB,$ and set $\mathscr{T}=\{x \in X: |x \rangle \langle x_kb| \in \Bw(\mathcal{A})\}.$ An easy calculation shows that $\sum_{i=1}^N x_i b_i \in \mathscr{T}$ for all $b_i \in \fB$ and $N \in \mathbb{N},$ and it follows from Lemma \ref{leftconv=rightconv} that $\mathscr{T}$ is WOT sequentially closed in $X$, so $\mathscr{T}=X.$ A similar argument show that $\{x \in X: |x \rangle \langle z| \in \Bw(\mathcal{A})\} = X$ for all $z \in X,$ and this proves the result.
\endproof

\begin{prop} \label{wcg}
Let $X$ be a right $C^*$-module over a $\bw$-algebra $\fB \subseteq B(\Hil).$ Then the following are equivalent:
\begin{enumerate}
\item[\rm{(1)}]  $X$ is $\swb$countably generated;
\item[\rm{(2)}]  $\mathbb{K}_\fB(X)$ has an element $T$ such that $\psi(T)>0$ for all nonzero WOT sequentially continuous positive functionals $\psi$ on $\Bw(\mathbb{K}_\fB(X))$;
\item[\rm{(3)}]  $\mathbb{K}_\fB(X)$ has a positive element $T$ with $\overline{T(X \otimes_\fB \Hil)} = X \otimes_\fB \Hil$;
\item[\rm{(4)}]  $\mathbb{K}_\fB(X)$ has an positive increasing sequential weak cai.
\end{enumerate}
If additionally $X$ is a $\bw$-module, these conditions imply that $\Bw(\mathbb{K}_\fB(X)) = \mathbb{B}_\fB(X).$
\end{prop}

\proof
The equivalence of (2), (3), and (4) follows from Lemma \ref{wsp lem}.

$(1) \implies (2)$ Suppose that $X$ is $\swb$countably generated by $\{x_i\}_{i=1}^\infty,$ and that these are scaled so that the series $\sum_{i=1}^\infty |x_i \rangle \langle x_i|$ converges in norm to a positive element $T$ in $\mathbb{K}_\fB(X).$ Let $\mathcal{A} \subseteq \K_\fB(\fX)$ be as in the proof of Lemma \ref{gen set}. By a calculation in the proof of \cite[Theorem 7.13]{BMP}, if $\varphi$ is a positive functional on $\Bw(\mathbb{K}_\fB(X)) \subseteq B(X \otimes_\fB \Hil)$ such that $\varphi(T)=0$, then $\varphi(a)=0$ for all $a \in \mathcal{A}.$ Let $\psi$ be a WOT sequentially continuous positive functional on $\Bw(\mathbb{K}_\fB(X))$ such that $\psi(T)=0.$ By the calculation just mentioned, $\text{Ker}(\psi)$ contains $\mathcal{A}$, and evidently $\text{Ker}(\psi)$ is sequentially WOT-closed. By Lemma \ref{gen set}, $\text{Ker}(\psi) = \Bw(\mathbb{K}_\fB(X)),$ so that $\psi = 0.$

$(4) \implies (1)$ Let $\{e_n\}_{n=1}^\infty$ be a weak cai for $\K_\fB(X).$ For each $n \in \N$, pick $x_i^n, y_i^n \in X$ for $i=1,...,m_n$ such that $\|\sum_{i=1}^{m_n} |x_i^n \rangle \langle y_i^n| - e_n \| < \frac 1n$ and $\|\sum_{i=1}^{m_n} |x_i^n \rangle \langle y_i^n | \| \leq 1.$ We claim that $f_n := \sum_{i=1}^{m_n} |x_i^n \rangle \langle y_i^n |$ is also a weak cai for $\K_\fB(X).$ To this end, let $K \in \K_\fB(X).$ Take two nonzero vectors $h,k \in X \otimes_\fB \Hil,$ let $\epsilon>0,$ and pick $N \in \N$ such that $\frac 1N < \frac \epsilon{2(\|K\|+1)\|h\|\|k\|}$ and $|\langle e_nKh,k \rangle - \langle Kh,k \rangle| < \frac \epsilon 2$ for all $n \geq N.$ Then for $n \geq N,$ $$| \langle f_nK h,k \rangle - \langle Kh,k \rangle | \leq | \langle f_n K h,k \rangle - \langle e_n Kh,k \rangle|+| \langle e_nK h,k \rangle - \langle K h,k \rangle| \leq \|f_n-e_n\| \|K\| \|h\| \|k\| + \frac \epsilon 2 < \epsilon.$$ Hence $f_nK \xrightarrow{\text{WOT}} K$, and so $\{f_n\}$ is a weak cai for $\K_\fB(X).$ By the final assertion in Proposition \ref{mapping spaces are bw}, $f_n |x \rangle \langle y|(z) \xrightarrow{WOT_\fB} |x \rangle \langle y| (z) = x \langle y|z \rangle$ for all $x,y,z \in X.$ But $f_n |x \rangle \langle y|(z) = \sum_{i=1}^{m_n} x_i^n \langle y_i^n|x \rangle \langle y|z \rangle,$ and so we have shown that every element in $X$ of the form $x \langle y|z \rangle$ is a WOT$_\fB$-limit of a sequence of elements from $\text{Span}\{x_i^n b: b \in \fB, n \in \N, i=1,...,m_n\}.$ Since the span of elements of the form $x \langle y|z \rangle$ is dense in $X$ (\cite[8.1.4 (2)]{BLM}), it follows that $X$ is WOT$_\fB$-generated by the countable set $\{x_i^n : n \in \N, i=1,...,m_n\}.$

For the last assertion, it follows directly from (4) that $I \in \Bw(\mathbb{K}_\fB(X)),$ and the assumption that $X$ is a $\bw$-module gives that $\mathbb{B}_\fB(X)$ is a $\bw$-algebra in $B(X \otimes_\fB \Hil),$ so that $\Bw(\mathbb{K}_\fB(X)) \subseteq \mathbb{B}_\fB(X)$. Since $\mathbb{K}_\fB(X)$ is an ideal in $\mathbb{B}_\fB(X)$, it follows that $\Bw(\mathbb{K}_\fB(X))$ is also an ideal in $\mathbb{B}_\fB(X)$, and so $\Bw(\mathbb{K}_\fB(X)) = \mathbb{B}_\fB(X).$
\endproof

\begin{prop} \label{cg adj ops}
Let $\fX$ be a $\swb$countably generated $\bw$-module over a $\bw$-algebra $\fB \subseteq B(\Hil).$ Then $B_\fB(\fX) = \mathbb{B}_\fB(\fX) = \Bw(\mathbb{K}_\fB(\fX)).$
\end{prop}

\proof
Suppose that $T \in B_\fB(\fX)$, and let $(e_n)$ be a sequence in $\mathbb{K}_\fB(\fX)$ with $e_n \nearrow I$ in $B(\fX \otimes_\fB \Hil).$ Viewing $T$ as an operator in $B(\fX \otimes_\fB \Hil),$ we have $Te_n \xrightarrow{\text{WOT}} T$ in $B(\fX \otimes_\fB \Hil).$ Since each $Te_n$ is adjointable, and $\mathbb{B}_\fB(\fX)$ is WOT sequentially closed in $B(\fX \otimes_\fB \Hil)$ by Proposition \ref{mapping spaces are bw}, we have proved that $T$ is adjointable. The last equality is the last assertion in Proposition \ref{wcg}.
\endproof

\begin{lemma} \label{inv sd}
If $\fX$ is a right $\bw$-module over a $\bw$-algebra $\fB \subseteq B(\Hil),$ then $\fX$ is selfdual as a $\fB$-module if and only if it is selfdual as a $\Bw(\langle X|X \rangle)$-module.
\end{lemma}

\proof
This follows directly from the general fact that for a $C^*$-module $X$ over $A$, $B_A(X,A)=B_J(X,J)$ for any ideal $J$ in $A$ containing $\langle X|X \rangle$ (see \cite[Lemma 8.5.2]{BLM}).
\endproof

\begin{thm} \label{bw=sd for cg}
Let $\fX$ be a $\swb$countably generated $C^*$-module over a $\bw$-algebra $\fB \subseteq B(\Hil).$ Then $\fX$ is a $\bw$-module over $\fB$ if and only if $\fX$ is selfdual.
\end{thm}

\proof
($\Rightarrow$) By Lemma \ref{inv sd}, we may assume without loss of generality that $\Bw(\langle X|X \rangle) = \fB.$ Let $\varphi \in B_\fB(\fX,\fB).$ Fix $x_0 \in X,$ and define $T:\fX \to \fX$ by $T(x)=x_0 \varphi(x)$ for $x \in \fX.$ It is easily checked that $T \in B_\fB(\fX)$, so by Proposition \ref{cg adj ops}, $T$ is adjointable, and by the easy direction of Proposition \ref{char of adj ops} (2), if $x_n \xrightarrow{WOT_\fB} x$ in $\fX$ and $y \in \fX,$ then
\[ \langle y|x_0 \rangle \varphi(x_n) = \langle y|T(x_n) \rangle \xrightarrow{\text{WOT}} \langle y|T(x) \rangle = \langle y|x_0 \rangle \varphi(x). \] Since $x_0$ was arbitrary, we have shown that for any $\zeta, \eta \in \Hil$ and $y,z \in \fX,$
\[ \langle \varphi(x_n) \zeta, \langle z|y \rangle \eta \rangle \to \langle \varphi(x) \zeta, \langle z|y \rangle \eta \rangle. \]
Hence $\varphi(x_n) \xrightarrow{\text{WOT}} \varphi(x),$ and so $\varphi = \langle y_0| \cdot \rangle$ for some $y_0 \in \fX$ by Proposition \ref{char of adj ops} (1).

$(\Leftarrow)$ Proposition \ref{selfdual implies bw}.
\endproof

\begin{cor} \label{cg is unital}
Let $\fB \subseteq B(\Hil)$ be a $\bw$-algebra considered as a $\bw$-module over itself. If $\fB$ is $\bw_\fB$-countably generated, then $\fB$ is unital.
\end{cor}

\proof
By Theorem \ref{bw=sd for cg}, $\fB$ is selfdual. Hence the identity map on $\fB$ is equal to $x \mapsto y^*x$ for some $y \in \fB$, so that $y$ is a unit for $\fB$.
\endproof


\begin{lemma} \label{dense perp0}
If $X$ is a WOT$_\fB$ sequentially dense subset of a $\bw$-module $\fX,$ then $X^\perp = (0).$
\end{lemma}

\proof
If $w \in X^\perp,$ then $\mathscr{S}=\{\xi \in \fX : \langle \xi|w \rangle = 0\}$ is WOT$_\fB$ sequentially closed and contains $X$, so $\mathscr{S}=\fX$. Hence $w \in \mathscr{S}$, so $w=0.$
\endproof

\begin{prop} \label{unique compl}
If $X$ is a (norm) countably generated $C^*$-module over a $\bw$-algebra $\fB \subseteq B(\Hil),$ then the $\bw$-module completion $\Bw(X)$ from Theorem \ref{completion} is the unique $\bw$-module containing $X$ as a WOT$_\fB$ sequentially dense submodule.
\end{prop}

\proof

(Cf.\ \cite{Hamana}, proof of uniqueness in Theorem 2.2) Let $\fY$ be another such $\bw$-module, and denote the $\fB$-valued inner product of $\fY$ by $(\cdot|\cdot).$ As in the proof of Proposition \ref{uniq}, define $U: \fY \to B(\Hil,X \otimes_\fB \Hil)$ by $U(\xi)^*(x) = (\xi|x)$ for $\xi \in \fY$ and $x \in X.$ It follows just as in that proof that $U$ is linear and $\fB$-linear and that $X \subseteq U(\fY) \subseteq \Bw(X).$ Also note that $U$ is bounded by the calculation $$\|U(\xi)\| = \|U(\xi)^*\| = \sup\{\|U(\xi)^*(x)\|: x \in \Ball(X)\} = \sup\{\|(\xi|x)\| : x \in \Ball(X)\} \leq \|\xi\|$$ for $\xi \in \fY.$ Now fix $\xi \in \fY,$ and consider the map $\fY \to \fB,$ $\eta \mapsto \langle U(\xi)|U(\eta) \rangle,$ which is easily seen to be in $B_\fB(\fY,\fB).$ Since $X$ is countably generated and WOT$_\fB$ sequentially dense in $\fY$, $\fY$ is a $\swb$countably generated $\bw$-module, and so is selfdual by Theorem \ref{bw=sd for cg}. Hence there is a $y_\xi \in \fY$ such that $$( y_\xi| \eta ) = \langle U(\xi)|U(\eta) \rangle \text{ for all } \eta \in \fY.$$ Define $T : \fY \to \fY$ by $T(\xi) = y_\xi,$ which is easily seen to be in $B_\fB(\fY) = \B_\fB(\fY).$ Consider $\Ker(\id_\fY - T) = \{y \in \fY: T(y)=y\}.$ This set contains $X$ and is WOT$_\fB$ sequentially closed since $\id_\fY-T$ is adjointable, so $\id_\fY = T.$ Thus $( \xi|\eta) = \langle U(\xi)|U(\eta) \rangle$ for all $\xi, \eta \in \fY.$

To prove that $U$ is a unitary between $\fY$ and $\Bw(X)$, it remains to show that $U(\fY)$ is WOT$_\fB$ sequentially closed in $\Bw(X).$ To this end, suppose that $(\xi_n)$ is a sequence in $\fY$ such that $U(\xi_n) \xrightarrow{WOT_\fB} \tau$ in $\Bw(X).$ By what we just proved, $$(\xi_n|\eta) = \langle U(\xi_n)|U(\eta) \rangle \xrightarrow{WOT} \langle \tau|U(\eta) \rangle \text{ for all } \eta \in \fY.$$ By Proposition \ref{bw char}, there exists a $\xi \in \fY$ such that $(\xi_n|\eta) \xrightarrow{WOT} (\xi|\eta)$ for all $\eta \in \fY.$ Thus $\langle \tau|U(\eta) \rangle = (\xi|\eta)=\langle U(\xi)|U(\eta) \rangle$; hence $\langle \tau-U(\xi)|U(\eta) \rangle = 0$ for all $\eta \in \fY.$ So $\tau-U(\xi) \in U(\fY)^\perp \subseteq X^\perp$. By Lemma \ref{dense perp0}, $X^\perp=(0),$ so $\tau = U(\xi).$
\endproof

In preparation for our analogue of Kasparov's stabilization theorem, we now present a direct sum construction for $\bw$-modules. Fix a $\bw$-algebra $\fB \subseteq B(\Hil),$ and let $\{\fX_n\}$ be a countable collection of $\bw$-algebras over $\fB.$ Define the \emph{direct sum $\bw$-module} to be the set $$ \oplus^w \fX_n := \{(x_n) \in \prod_n \fX_n : \sum_n \langle x_n|x_n \rangle \text{ is WOT-convergent in } \fB\},$$ with the inner product $\langle (x_n)|(y_n) \rangle = \sum_n \langle x_n| y_n \rangle$ and obvious $\fB$-module action. The proof below that this is a $\bw$-module follows \cite[8.5.26]{BLM} pretty closely.

\begin{lemma} \label{sum is bw mod}
The space $\fX := \oplus^w \fX_n$ defined above is a $\bw$-module over $\fB.$
\end{lemma}

\proof
It follows from the operator inequality $(y+i^k x)^*(y + i^k x) \leq 2(x^*x+y^*y)$ and the polarization identity $x^*y = \frac 14 \sum_{k=0}^3 i^k (y+i^kx)^*(y+i^k x)$ that $\langle (x_n)|(y_n) \rangle := \sum_n \langle x_n| y_n \rangle$ does indeed define a $\fB$-valued inner product on $\fX.$ It is easy to check that $\langle \cdot|\cdot \rangle$ satisfies all the axioms of a $C^*$-module inner product.

Define $\Kil_n := \fX_n \otimes_\fB \Hil$ and $\Kil := \oplus^2_n \Kil_n,$ and let $P_n: \Kil \to \Kil_n$ be the canonical projection. Since each $\fX_n$ is WOT sequentially closed in $B(\Hil,\Kil_n)$, it follows immediately that the space $$W:=\{T \in B(\Hil, \Kil) : P_nT \in \fX_n \text{ for all }n\}$$ is a WOT sequentially closed TRO. By Theorem \ref{pics}, $W$ is a $\bw$-module over $\Bw(W^*W) \subseteq \fB.$ Hence $W$ is a $\bw$-module over $\fB$. Define a $\fB$-module map $U: \fX \to W$ by sending $(x_n) \in \fX$ to the SOT-convergent sum $\sum_n P_n^* x_n$ (indeed, for $\zeta \in \Hil$ and $N,M \in \N$ with $N \geq M,$ a short calculation gives $\|\sum_{n=M}^N P_n^* x_n(\zeta)\|^2 = \sum_{n=M}^N \langle \zeta, \langle x_n|x_n \rangle \zeta \rangle,$ and so by Cauchy's convergence test, the series $\sum_n P_n^* x_n (\zeta)$ converges). To check surjectivity of $U$, note that $\{P_n^* P_n\}_{n=1}^\infty$ is a family of mutually orthogonal projections in $B(\Kil)$ with $\sum_n P_n^* P_n = I_\Kil.$ If $T \in W,$ then $\sum_{n=1}^N \langle P_n T| P_n T \rangle = T^*(\sum_{n=1}^N P_n^* P_n) T \leq T^*T.$ Thus $(P_n T) \in \fX$, and $U((P_nT)) = \sum_n P_n^* P_n T = T.$ Finally, for $x,y \in \fX,$ the formula $$U(x)^*U(y) = \langle x|y \rangle$$ is an easy exercise (first checking this when both $x$ and $y$ are ``finitely supported," then extending via WOT-limits to the general case). 

So we have established the existence of a surjective inner-product-preserving $\fB$-module map $U: \fX \to W,$ where $W$ is a $\bw$-module over $\fB.$ It follows immediately that $\fX$ is complete, and a straightforward application of Proposition \ref{bw char} shows that $\fX$ is also a $\bw$-module over $\fB.$
\endproof


Letting $\fX_n = \fB$ for all $n$, we obtain a $\bw$-module over $\fB$ which we denote as $C^w(\fB)$.

\begin{cor} \label{col sd}
If $\fB$ is a unital $\bw$-algebra, then the $\bw$-module $C^w(\fB)$ is selfdual.
\end{cor}

\proof
Since $\fB$ is unital, $C^w(\fB)$ is $\swb$countably generated.
\endproof

Early in this investigation, David Blecher proved an interesting generalization of Corollary \ref{col sd}, and we thank him for allowing it to be included here.

For a cardinal number $I$, a $C^*$-algebra $B \subseteq B(\Hil)$ is said to be \emph{$I$-additively weakly closed} if whenever $\sum_{k \in I} x_k^* x_k$ is bounded in $B(\Hil)$ for a collection $\{x_k\}_{k \in I}$ in $B,$ then the WOT-limit of this sum is an element in $B$. For an $I$-additively weakly closed $B \subseteq B(\Hil),$ define $$C_I^w(B) = \{(x_k) \in \prod_{k \in I} B : \sum_k x_k^* x_k \text{ is WOT-convergent in } B \}.$$ One may then show as in the first part of the proof of Lemma \ref{sum is bw mod} that $\langle (x_k)|(y_k) \rangle := \sum_{k \in I} x_k^* y_k$ defines a $B$-valued inner product on $C_I^w(B)$. It is easy to argue that this satisfies all the axioms of a $C^*$-module inner product. Completeness of $C_I^w(B)$ follows as in the second to last paragraph of \cite[1.2.26]{BLM}, since $C_I^w(B)$ clearly coincides with the (underlying Banach space of the) operator space of the same notation there.

To see that Proposition \ref{dpb} below generalizes Corollary \ref{col sd}, note that an easy ``telescoping series" argument shows that $B$ is $\mathbb{N}$-additively weakly closed if and only if $B$ is a Borel $*$-algebra (that is, closed under weak limits in $B(\Hil)$ of bounded monotone sequences of selfadjoint elements). Hence every $\bw$-algebra is $\mathbb{N}$-additively weakly closed.

To set some notation for the following lemma and proposition, let $B \subseteq B(\Hil)$ be a nondegenerate $I$-additively weakly closed $C^*$-algebra. For each $j \in I,$ denote by $\varepsilon_j : \Hil \to \Hil^{(I)}$ the canonical inclusion into the $j^{th}$ summand, and by $P_j : \Hil^{(I)} \to \Hil$ the canonical projection from the $j^{th}$ summand (so $\varepsilon_j^* = P_j$). For $b \in B$ and $j \in I,$ denote by $e_jb$ the element in $C_I^w(B)$ with $b$ in the $j^{th}$ slot and 0's elsewhere.

\begin{lemma} \label{addcl lemma}
If $B \subseteq B(\Hil)$ is an $I$-additively weakly closed $C^*$-algebra, then $C_I^w(B) \otimes_B \Hil \cong \Hil^{(I)}$ via a unitary $U: \Hil^{(I)} \to C_I^w(B) \otimes_B \Hil$ such that $$U(\varepsilon_j(b\zeta)) = e_jb \otimes \zeta \text{ for all }b \in B,\ j \in I, \text{ and } \zeta \in \Hil$$ and $$U^*((b_i) \otimes \zeta) = (b_i \zeta) \text{ for all } (b_i) \in C_I^w(B) \text{ and } \zeta \in \Hil.$$
\end{lemma}

\proof
By Cohen's factorization theorem (\cite[A.6.2]{BLM}), every element in $\Hil$ can be expressed in the form $b \zeta$ for some $b \in B$ and $\zeta \in \Hil.$ Using this, define a map $U_0: \mathcal{F} \to C_I^w(B) \otimes_B \Hil$ on the dense subspace $\mathcal{F}$ of finitely supported columns in $\Hil^{(I)}$ by $$U_0(\sum_{j \in F} \varepsilon_j(b_j \zeta_j)) = \sum_{j \in F} e_jb_j \otimes \zeta_j$$ for $(b_j \zeta_j) \in \mathcal{F}$ supported on a finite subset $F \subseteq I.$ To see that this is well-defined, suppose that $b,b' \in B$ and $\zeta, \zeta' \in \Hil$ with $b \zeta = b' \zeta'$. Then for any $(c_i) \otimes \eta \in C_I^w(B) \otimes_B \Hil,$ $$\langle e_jb \otimes \zeta - e_j b' \otimes \zeta', (c_i) \otimes \eta \rangle = \langle \zeta, b^*c_j \eta \rangle - \langle \zeta', (b')^*c_j \eta \rangle = \langle b \zeta - b' \zeta', c_j \eta \rangle = 0.$$ By totality of the simple tensors in $C_I^w(B) \otimes_B \Hil,$ $e_j b \otimes \zeta - e_j b' \otimes \zeta' = 0.$ It follows that $U_0$ is well-defined. A direct calculation shows that $U_0$ is isometric, hence extends to an isometry $U: \Hil^{(I)} \to C_I^w(B) \otimes_B \Hil.$ To see that $U$ is surjective, let $(b_i) \in C_I^w(B)$, $\zeta \in \Hil,$ take $F \subseteq I$ to be finite, and denote by $(b_i)_F$ the ``restriction" of $(b_i)$ to $F$. Then
\begin{equation*}
\begin{split}
\|(b_i) \otimes \zeta - \sum_{i \in F} e_ib_i \otimes \zeta \|^2 &= \langle \zeta, \langle (b_i)-(b_i)_F|(b_i)-(b_i)_F \rangle \zeta \rangle \\
	&= \langle \zeta, (\sum_{i \in I} b_i^* b_i - \sum_{i \in F} b_i^* b_i) \zeta \rangle.
\end{split}
\end{equation*}
If we interpret $(\sum_{i \in I} b_i^* b_i - \sum_{i \in F} b_i^* b_i)$ as a net indexed by the collection of finite subsets $F$ of $I$, the last displayed quantity converges to 0.
So $$U(\sum_{i \in F} \varepsilon_i(b_i \zeta)) = \sum_{i \in F} e_i b_i \otimes \zeta \longrightarrow (b_i) \otimes \zeta \text{ in norm},$$ where $\sum_{i \in F} e_ib_i \otimes \zeta$ is considered to be a net indexed by the collection of finite subsets $F$ of $I$. Since the set of simple tensors in $C_I^w(B) \otimes_\fB \Hil$ spans a dense subset, it follows that $U$ is surjective. The first displayed equation in the statement is obvious from the first displayed equation in this proof. For the second, we need to show $\langle (b_i) \otimes \zeta, U((\zeta_i)) \rangle = \langle (b_i \zeta), (\zeta_i) \rangle$ for all $(b_i) \in C_I^w(B),\ \zeta \in \Hil,$ and $(\zeta_i) \in \Hil^{(I)},$ which we leave as an exercise by first checking for finitely supported $(\zeta_i).$
\endproof

\begin{prop}[David Blecher] \label{dpb}
If $B \subseteq B(\Hil)$ is a unital and $I$-additively weakly closed $C^*$-algebra, then the $C^*$-module $C_I^w(B)$ is selfdual.
\end{prop}

\proof
We first fix some notation. Following \cite[1.2.26]{BLM}, denote by $\mathbb M_I(B(\Hil))$ the space of $I \times I$ matrices over $B(\Hil)$ whose finite submatrices have uniformly bounded norm, and equip this space with the norm $$\|u\| = \sup\{ \|u_F\| : u_F \text{ is a finite submatrix of } u\}.$$ It is well-known (see e.g.\ the section in \cite{BLM} just mentioned) that this is a Banach space that is canonically isometrically isomorphic to $B(\Hil^{(I)})$. Denote by $\mathbb M_I(B)$ the subspace of $\mathbb M_I(B(\Hil))$ consisting of matrices with entries in $B$.

Fixing an index $j_0 \in I$, there is a canonical isometric embedding of $C_I^w(B)$ onto the subspace of $\mathbb M_I(B)$ consisting of matrices supported on the $j_0^{\mathrm{th}}$ column (we omit the routine details of this), and a canonical embedding of $B$ onto the subspace of matrices in $\mathbb M_I(B)$ supported on the $(j_0,j_0)$-entry. Write
\begin{align*}
\rho: C_I^w(B) \hookrightarrow \mathbb M_I(B) \\
\sigma: B \hookrightarrow \mathbb M_I(B)
\end{align*} for these embeddings.

We show that there is also a canonical embedding $$\pi: B_B(C_I^w(B)) \hookrightarrow \mathbb M_I(B).$$ Indeed, by Proposition \ref{bigdiagram}, Lemma \ref{addcl lemma}, and \cite{BLM} (1.19), we have the following canonical embedding and isomorphisms:
$$B_B(C_I^w(B)) \hookrightarrow B(C_I^w(B) \otimes_B \Hil) \cong B(\Hil^{(I)}) \cong \mathbb{M}_I(B(\Hil)).$$
Using the unitary $U$ from Lemma \ref{addcl lemma}, we have $$P_i \circ U^*TU \circ \varepsilon_j (\zeta) = P_i(U^*T(e_j \otimes \zeta)) = P_i(U^*(T(e_j) \otimes \zeta)) = (Te_j)_i(\zeta) = \langle e_i|T(e_j) \rangle (\zeta)$$ for all $i,j \in I$ and $\zeta \in \Hil.$ That is, under the embedding and isomorphisms just mentioned, $T \in B_B(C_I^w(B))$ corresponds to the matrix $[T_{ij}] \in \mathbb{M}_I(B(\Hil))$ with $T_{ij} = \langle e_i|T(e_j) \rangle \in B.$

It is straightforward (using the definitions of $\rho$, $\pi$, and $\sigma$ as composite maps involving the unitary $U$ from Lemma \ref{addcl lemma}) to show that for $x,y \in C_I^w(B)$ and $T \in B_B(C_I^w(B))$, we have 
\begin{align*}
\rho(Tx) = \pi(T) \rho(x), \\
\sigma(\langle x|y \rangle) = \rho(x)^* \rho(y).
\end{align*}
Note also that $\pi(T)^* \rho(y)$ is a matrix in $\mathbb M_I(B)$ supported on the $j_0^{\mathrm{th}}$ column, so $\pi(T)^* \rho(y) = \rho(z)$ for some $z \in C_I^w(B)$.
Hence \[ \sigma(\langle Tx|y \rangle) = \rho(Tx)^* \rho(y) = (\pi(T) \rho(x))^* \rho(y) = \rho(x)^* \pi(T)^* \rho(y) = \rho(x)^* \rho(z) = \sigma(\langle x|z \rangle). \]
So $\langle Tx|y \rangle = \langle x|z \rangle$, and this is enough to prove that $T$ is adjointable.

Thus $B_B(C_I^w(B)) = \mathbb B_B(C_I^w(B))$. To prove selfduality, let $\tau \in B_B(C_I^w(B),B)$ and fix an index $k \in I$. Define $T \in B_B(C_I^w(B))$ by $T(x) = e_k \tau(x)$. Then \[ \tau(x) = \langle e_k | e_k \tau(x) \rangle = \langle e_k | T(x) \rangle = \langle T^*(e_k) | x \rangle \] for all $x \in C_I^w(B)$, so that $\tau = \langle T^*(e_k) | \cdot \rangle$.
\endproof

The following lemma is a $\bw$-analogue of Lemma 2.34 from \cite{RW} or Proposition 3.8 in \cite{Lan}. Note that the simple proof presented in these books does not seem to work in our setting, since it is unclear how to extend an isometry from a WOT$_\fB$ sequentially dense subspace to the whole space. In the proof below, we write $\Bw(S)$ to denote the WOT$_\fB$ sequential closure of a subset $S$ of a $\bw$-module over $\fB$. (Recall from Note \ref{dist notation} that for a sequence in a $\bw$-module $\fX$ over $\fB,$ WOT$_\fB$-convergence coincides with WOT-convergence in $B(\Hil, \fX \otimes_\fB \Hil),$ so this notation does not clash with our previous meaning of $\Bw(\cdot)$ as WOT sequential closure of subsets of $B(\Hil).$)

\begin{lemma} \label{unitary equiv}
Let $\fX, \fY$ be $\bw$-modules over $\fB \subseteq B(\Hil).$ If $T$ is an operator in $\B_\fB(\fX,\fY)$ such that $T(\fX)$ is WOT$_\fB$ sequentially dense in $\fY$ and $T^*(\fY)$ is WOT$_\fB$ sequentially dense in $\fX,$ then $\fX$ and $\fY$ are unitarily equivalent.
\end{lemma}

\proof
Consider $T$ as an element in the $\bw$-algebra $\B_\fB(\fX \oplus \fY) \cong \left[ \begin{matrix} \B_\fB(\fX) & \B_\fB(\fY,\fX) \\ \B_\fB(\fX,\fY) & \B_\fB(\fY) \end{matrix} \right] \subseteq B((\fX \otimes_\fB \Hil) \oplus(\fY \otimes_\fB \Hil)),$ and take the polar decomposition $$\left[ \begin{matrix} 0 & 0 \\ T & 0 \end{matrix} \right] = U \left[ \begin{matrix} |T| & 0 \\ 0 & 0 \end{matrix} \right] = \left[ \begin{matrix} U_{11} & U_{12} \\ U_{21} & U_{22} \end{matrix} \right] \left[ \begin{matrix} |T| & 0 \\ 0 & 0 \end{matrix} \right].$$ By Proposition \ref{polar decomp}, $U \in \B_\fB(\fX \oplus \fY)$. We see that $U_{11}|T|=0$, and the formula $U^* T = |T|$ shows that $U_{22}^* T=0.$ Consider the set $\mathscr{T} = \{y \in \fY : T^*(y) \in \Bw(|T|\fX)\}$, where $\Bw(S)$ denotes the WOT$_\fB$ sequential closure of a subset $S \subseteq \fX$. Since $T^*$ is adjointable, it is WOT$_\fB$ sequentially continuous, so $\mathscr{T}$ is a WOT$_\fB$ sequentially closed subset of $\fY$ containing $T(\fX)$. Hence $\mathscr{T}=\fY,$ i.e.\ $T^*(\fY) \subseteq \Bw(|T|\fX),$ so $\fX = \Bw(T^*(\fY)) \subseteq \Bw(|T|\fX)$ (using the notation mentioned above the statement of the lemma). Since $U_{11}$ is WOT$_\fB$ sequentially continuous, its kernel in $\fX$ is WOT$_\fB$ sequentially closed, and since $\Ker(U_{11})$ contains the WOT$_\fB$ sequentially dense set $|T|(\fX),$ we have $U_{11}=0.$ A similar but shorter argument shows that $U_{22}=0$ as well.

Since $U$ is a partial isometry, it follows now that $U_{21}$ is as well. The relation $U_{21}=U_{21}U_{21}^*U_{21}$ implies that $U_{21}(\fX)$ is WOT$_\fB$ sequentially closed in $\fY$. (Indeed, suppose $U_{21}(x_n) \xrightarrow{WOT} y$ in $\fY$. Then $U_{21}(x_n)=U_{21}U_{21}^*U_{21}(x_n) \xrightarrow{WOT} U_{21}U_{21}^*(y)$ since $U_{21}$ and $U_{21}^*$ are adjointable, hence WOT$_\fB$ sequentially continuous. Thus $y = U_{21}U_{21}^*(y) \in U_{21}(\fX).$) Since $T=U_{21}|T|$, $U_{21}(\fX)$ contains the WOT$_\fB$ sequentially dense set $T(\fX)$, and so $U_{21}$ is surjective. Similarly, $U_{21}^*$ is a partial isometry with $U_{21}^*(\fY) = \fX,$ and it follows that $U_{21}: \fX \to \fY$ is a unitary.
\endproof

It is quite surprising that (given Lemma \ref{unitary equiv}) the obvious $\bw$-analogue of Kasparov's stabilization theorem now follows from only a slight modification of the proof presented in \cite[Theorem 6.2]{Lan} and \cite[Theorem 5.49]{RW} for the original stabilization theorem.

\begin{thm} \label{stabilization}
If $\fB \subseteq B(\Hil)$ is a $\bw$-algebra and $\fX$ is a $\swb$countably generated $\bw$-module over $\fB,$ then $\fX \oplus C^w(\fB) \cong C^w(\fB)$ unitarily.
\end{thm}

\proof
Using the second comment in the paragraph under Definition \ref{ssmod def} to make sense of the reduction to the unital case, apply the argument in \cite[Theorem 6.2]{Lan} or \cite[Theorem 5.49]{RW}, changing ``generating set" to ``$\swb$generating set," ``dense" to ``WOT$_\fB$ sequentially dense," and the $C^*$-module direct sum of countably many copies of $\fB$ to $C^w(\fB).$ Finish off the argument by invoking Lemma \ref{unitary equiv}.
\endproof

We now discuss more generally orthogonally complemented submodules of $\bw$-modules, and then make a connection between a $\bw$-analogue of the $C^*$-module theory of quasibases and orthogonally complemented submodules of the $\bw$-module $C^w(\fB).$

\begin{defn}
\begin{itemize}
\item A closed submodule $X$ of a $C^*$-module $Y$ over a $C^*$-algebra $A$ is said to be \emph{orthogonally complemented in $Y$} if there is another closed submodule $W$ in $Y$ such that $W+X = Y$ and $\langle w|x \rangle=0$ for all $w \in W$ and $x \in X.$ It is easy to see that this happens exactly when $X$ is the range of a projection $P \in \B_A(Y).$ (For one direction of this, check that if $X$ is orthogonally complemented in $Y$ with $W$ as above, then each element in $Y$ has a unique representation as a sum $x+w$ with $x \in X,$ $w \in W.$ Then show that the map $P: Y \to Y$ defined $P(w+x)=x$ for $w \in W$ and $x \in X$, satisfies $\langle P(x+w)|x'+w' \rangle = \langle x+w|P(x'+w') \rangle$ for $w,w' \in W,$ $x,x' \in X.$)

\item A closed submodule $\fX$ of a $\bw$-module $\fY$ over $\fB$ will be called a \emph{$\bw$-submodule of $\fY$} if $\fX$ is a $\bw$-module with the $C^*$-module structure it inherits from $\fY$, and if $\fX$ satisfies the following additional condition: whenever $(x_n)$ is a sequence in $\fX$ and $x \in \fX$ such that $\langle x_n|z \rangle \xrightarrow{WOT} \langle x|z \rangle$ for all $z \in \fX$, then $\langle x_n|y \rangle \xrightarrow{WOT} \langle x|y \rangle$ for all $y \in \fY$ (in other words, if a sequence converges WOT$_\fB$ in $\fX$, then it converges WOT$_\fB$ in $\fY$ to the same limit).
\end{itemize}
\end{defn}

\begin{prop} \label{orth compl}
Let $\fX$ be a closed submodule of a $\bw$-module $\fY$ over a $\bw$-algebra $\fB \subseteq B(\Hil).$ Consider the following conditions:
\begin{enumerate}
\item[\rm{(1)}]  $\fX$ is orthogonally complemented in $\fY$;

\item[\rm{(2)}]  $\fX$ is a $\bw$-submodule of $\fY$;

\item[\rm{(3)}]  $\fX$ is a $\bw$-module with the inner product and module structure inherited from $\fY$, and $\fX$ is WOT$_\fB$ sequentially closed in $\fY.$

We have $(1) \implies (2) \iff (3).$ If $\fX$ satisfies $(2)$ and $\Bw(\K_\fB(\fX))=\B_\fB(\fX),$ then $(1)$ holds.
\end{enumerate}
\end{prop}

\proof
$(1) \implies (2).$ Let $P \in \B_\fB(\fY)$ be a projection with range $\fX.$ Suppose $(x_n)$ is a sequence in $\fX$ such that $\langle x_n|w \rangle$ is WOT-convergent for all $w \in \fX.$ Since $\fY$ is a $\bw$-module and $\langle x_n|y \rangle = \langle x_n|Py \rangle$ for all $y \in \fY$, it follows from Proposition \ref{bw char} that there is an $x \in \fY$ such that $\langle x_n|y \rangle \xrightarrow{WOT} \langle x|y \rangle$ for all $y \in \fY.$ But then $\langle x_n|y \rangle = \langle x_n|Py \rangle \xrightarrow{WOT} \langle x|Py \rangle = \langle Px|y \rangle$ for all $y \in \fY,$ so $x=Px \in \fX.$ This proves that $\fX$ is a $\bw$-module. To see that $\fX$ is a $\bw$-submodule of $\fY,$ suppose that $(x_n),x \in \fX$ with $\langle x_n | w \rangle \xrightarrow{WOT} \langle x|w \rangle$ for all $w \in \fX.$ Then $\langle x_n|y \rangle = \langle x_n|Py \rangle \xrightarrow{WOT} \langle x|Py \rangle = \langle x|y \rangle$ for all $y \in \fY.$

$(2) \implies (3).$ By definition, a $\bw$-submodule is a $\bw$-module with the inherited structure. To show that $\fX$ is WOT$_\fB$ sequentially closed in $\fY,$ suppose that $(x_n)$ is a sequence in $\fX$ converging WOT$_\fB$ to $y$ in $\fY,$ i.e.\ $\langle x_n|w \rangle \xrightarrow{WOT} \langle y|w \rangle$ for all $w \in \fY.$ In particular, $\langle x_n|z \rangle$ is WOT-convergent for all $z \in \fX,$ so by Proposition \ref{bw char}, there is an $x \in \fX$ such that $\langle x_n|z \rangle \xrightarrow{WOT} \langle x|z \rangle$ for all $z \in \fX.$ By the ``additional condition" in the definition of $\bw$-submodule, $\langle x_n|w \rangle \xrightarrow{WOT} \langle x|w \rangle$ for all $w \in \fY.$ Hence $y = x \in \fX.$

$(3) \implies (2).$ Assume (3). Note that we can canonically identify $\fX \otimes_\fB \Hil$ with a closed subspace of $\fY \otimes_\fB \Hil.$ Indeed, the canonical inclusion of $\{\sum_{i=1}^n x_i \otimes \zeta_i \in \fX \otimes_\fB \Hil: x_i \in \fX, \zeta_i \in \Hil\}$ into $\fY \otimes_\fB \Hil$ is isometric, hence extends to an isometry from $\fX \otimes_\fB \Hil$ into $\fY \otimes_\fB \Hil.$ To see that $\fX$ is a $\bw$-module with the inherited $C^*$-module structure, suppose that $(x_n)$ is a sequence in $\fX$ viewed in $B(\Hil, \fX \otimes_\fB \Hil)$ with $x_n \xrightarrow{WOT} T$ in $B(\Hil, \fX \otimes_\fB \Hil)$. Note that there is a canonical WOT-continuous embedding of $B(\Hil, \fX \otimes_\fB \Hil)$ into $B(\Hil, \fY \otimes_\fB \Hil)$ making the following diagram commute:
\[
\begin{tikzcd}
B(\Hil, \fX \otimes_\fB \Hil) \arrow{r} & B(\Hil, \fY \otimes_\fB \Hil) \\ \fX \arrow{u} \arrow{r} & \fY \arrow{u}
\end{tikzcd}
\]
So $x_n \xrightarrow{WOT} T$ in $B(\Hil, \fY \otimes_\fB \Hil),$ and since $\fY$ is WOT sequentially closed in the latter, $T \in \fY$ and $x_n \xrightarrow{WOT_\fB} T$ in $\fY.$ By the assumption that $\fX$ is WOT$_\fB$ sequentially closed in $\fY,$ $T \in \fX.$ Hence $\fX$ is WOT sequentially closed in $B(\Hil,\fX \otimes_\fB \Hil)$, and so by definition, $\fX$ is a $\bw$-module.

Now suppose that $(x_n),x$ are in $\fX$ and $\langle x_n|z \rangle \xrightarrow{WOT} \langle x|z \rangle$ for all $z \in \fX.$ Fixing $\zeta, \eta \in \Hil,$ we have $$\langle x_n \otimes \zeta, z \otimes \eta \rangle = \langle \zeta, \langle x_n|z \rangle \eta \rangle \to \langle \zeta, \langle x_n|z \rangle \eta \rangle = \langle x_n \otimes \zeta, z \otimes \eta \rangle$$ for all $z \in \fX.$ Take $y \in \fY,$ and let $\epsilon>0.$ Denote by $P$ the projection in $B(\fY \otimes_\fB \Hil)$ with range $\fX \otimes_\fB \Hil.$ By the principle of uniform boundedness, there is a $K>0$ such that $\|x_n\| \leq K$ for all $n$ and $\|x\| \leq K.$ Pick $\sum_{i=1}^k z_i \otimes \zeta_i \in \fX \otimes_\fB \Hil$ with $$\|P(y \otimes \eta) - \sum_{i=1}^k z_i \otimes \zeta_i\| < \frac{\epsilon}{3K(\|\zeta\|+1)},$$ and pick $N \in \N$ with $$|\langle x_n \otimes \zeta, \sum_{i=1}^k z_i \otimes \zeta_i \rangle -\langle x \otimes \zeta, \sum_{i=1}^k z_i \otimes \zeta_i \rangle| < \frac{\epsilon}3$$ for all $n \geq N.$ A triangle inequality argument then gives $$|\langle \zeta, \langle x_n|y \rangle \eta \rangle - \langle \zeta, \langle x|y \rangle \eta \rangle | = |\langle x_n \otimes \zeta, P(y \otimes \eta) \rangle - \langle x \otimes \zeta, P(y \otimes \eta) \rangle  < \epsilon$$ for all $n \geq N.$ Since $\zeta, \eta \in \Hil$ were arbitrary, we have shown that $\langle x_n |y \rangle \xrightarrow{WOT} \langle x|y \rangle$ for all $y \in \fY.$

Now we prove the final claim in the statement of the proposition. Suppose that $\fX$ is a $\bw$-submodule of $\fY$ and that $\Bw(\K_\fB(\fX))=\B_\fB(\fX).$ By definition of $\bw$-submodule, the inclusion $\iota: \fX \hookrightarrow \fY$ is WOT$_\fB$ sequentially continuous, so by Proposition \ref{char of adj ops} (2), $\iota$ is adjointable. Since $\langle \iota^* \iota x|x' \rangle = \langle \iota x| \iota x' \rangle = \langle x|x' \rangle$ for all $x,x' \in \fX,$ we have $\iota^* \iota = \rm{id}_\fX.$ It follows that $P=\iota \iota^* \in \B_\fB(\fY)$ is a projection with range $\fX.$
\endproof

\begin{note}
(We thank David Blecher for pointing this out.) To show that (2)/(3) does not imply (1) in Proposition \ref{orth compl} in general, let $\fB \subseteq B(\Hil)$ be a nonunital $\bw$-algebra, and take $\fX$ to be $\fB$ and $\fY$ to be the unitization $\fB^1 \subseteq B(\Hil),$ where we view these both as $\bw$-modules over $\fB^1.$ Clearly $\fX$ satisfies (3), but $\fX$ is not orthogonally complemented in $\fY$ since it is a proper subset and $\{y \in \fY : \langle y|x \rangle=0 \text{ for all } x \in \fX\} = \{c+\mu I_\Hil \in \fB^1 : b^*c + \mu b^* = 0 \text{ for all } b \in \fB\} = (0)$ by an approximate identity argument.
\end{note}

We now define a $\bw$-analogue of ``quasibases" for $C^*$-modules (see \cite{BLM} 8.2.5 and relevant notes in 8.7).

\begin{defn}
For a $\bw$-module $\fX$ over a $\bw$-algebra $\fB \subseteq B(\Hil)$, a countable subset $\{x_k\}$ of $\fX$ is a called a \emph{weak quasibasis} for $\fX$ if for any $x \in \fX,$ the sequence of finite sums $\sum_{k=1}^n x_k \langle x_k|x \rangle$ WOT$_\fB$-converges to $x.$ In other words, $\{x_k\}$ is a weak quasibasis iff $\sum_{k=1}^n |x_k \rangle \langle x_k | \nearrow I$ in $B(\fX \otimes_\fB \Hil).$
\end{defn}

\begin{remark} Quasibases are also called ``frames" (it appears that ``quasibasis" is the older term for these and ``frame" is the term most commonly employed in recent literature). Frank and Larson in \cite{FrankLar} initiated a systematic study of quasibases/frames for Hilbert $C^*$-modules, and what we have called ``weak quasibases" are essentially equivalent to ``non-standard normalized tight frames" in the terminology of their paper (see \cite[Definition 2.1]{FrankLar}). We also remark that Frank and Larson followed a similar approach to ours in using Kasparov's stabilization theorem to deduce the existence of quasibases/frames. (We thank the referee for drawing our attention to these points.)
\end{remark}

\begin{prop} \label{oc and qb}
Let $\fB \subseteq B(\Hil)$ be a $\bw$-algebra, and let $\fX$ be a right Banach module over $\fB.$

If $\fX$ is a $\bw$-module over $\fB$ with a weak quasibasis, then $\fX$ is isometrically $\fB$-isomorphic to an orthogonally complemented submodule of $C^w(\fB).$

Conversely, if $\fX$ is isometrically $\fB$-isomorphic to an orthogonally complemented submodule of a $\bw$-module $\fY$ over $\fB$ such that $\fY$ has a weak quasibasis, then $\fX$ is a $\bw$-module over $\fB$ with the canonically induced inner product, and $\fX$ has a weak quasibasis.
\end{prop}

\proof
For the first statement, if $\fX$ is a $\bw$-module over $\fB$ with a weak quasibasis, then clearly $\fX$ is $\swb$countably generated, and the result now follows from Theorem \ref{stabilization}.

For the converse, note first that if $\fX$ is isometrically $\fB$-isomorphic to any $\bw$-module $\fX_0$ over $\fB$ via an isometric $\fB$-isomorphism $U: \fX \to \fX_0,$ then defining $\langle x|y \rangle := \langle Ux|Uy \rangle$ for $x,y \in \fX$ makes $U$ a unitary between $C^*$-modules and makes $\fX$ a $\bw$-module over $\fB$ (this can be seen either by applying Proposition \ref{bw char} or invoking Definition \ref{ssmod def}, noting that $U$ induces a canonical unitary $\fX \otimes_\fB \Hil \cong \fX_0 \otimes_\fB \Hil$). Since every orthogonally complemented submodule of $\fY$ is a $\bw$-module over $\fB$ by Proposition \ref{orth compl}, and since unitaries between $\bw$-modules preserve weak quasibases, we may assume without loss of generality that $\fX$ is actually an orthogonally complemented submodule of $\fY.$ In that case, let $P \in \B_\fB(\fY)$ be a projection with range $\fX,$ and let $\{e_k\}$ be a weak quasibasis for $\fY.$ Then for any $x \in \fX,$ $$\sum_{k=1}^n P(e_k) \langle P(e_k)|x \rangle = \sum_{k=1}^n P(e_k) \langle e_k|x \rangle = P(\sum_{k=1}^n e_k \langle e_k|x \rangle) \xrightarrow{WOT_\fB} P(x)=x,$$ so that $\{P(e_k)\}$ is a weak quasibasis for $\fX.$
\endproof

We close by coalescing some of the main results of this section in the case of a unital coefficient $\bw$-algebra:

\begin{thm}
Let $\fB \subseteq B(\Hil)$ be a unital $\bw$-algebra, and let $\fX$ be a Banach module over $\fB.$ The following are equivalent:
\begin{enumerate}
\item[\rm{(1)}]  $\fX$ is a $\swb$countably generated $\bw$-module over $\fB$;

\item[\rm{(2)}]   $\fX$ is a $\swb$countably generated selfdual $C^*$-module over $\fB$;

\item[\rm{(3)}]   $\fX$ is a $\bw$-module with a weak quasibasis;

\item[\rm{(4)}]   $\fX$ is isometrically $\fB$-isomorphic to an orthogonally complemented submodule of $C^w(\fB)$;

\item[\rm{(5)}]   $\fX \oplus C^w(\fB) \cong C^w(\fB).$
\end{enumerate}
\end{thm}

\proof
$(1) \iff (2).$ Theorem \ref{bw=sd for cg}.

$(1) \implies (5).$ Theorem \ref{stabilization}.

$(5) \implies (4).$ Easy.

$(4) \implies (3).$ Proposition \ref{oc and qb}, noting that if $\fB$ is unital, then $C^w(\fB)$ has a canonical weak quasibasis.

$(3) \implies (1).$ Easy.
\endproof

\end{document}